%

\magnification=\magstep1
\def\forces{\parallel\!\!\! -}


\def\hexnumber#1{\ifcase#1 0\or1\or2\or3\or4\or5\or6\or7\or8\or9\or
	A\or B\or C\or D\or E\or F\fi }

\font\teneuf=eufm10
\font\seveneuf=eufm7
\font\fiveeuf=eufm5
\newfam\euffam
\textfont\euffam=\teneuf
\scriptfont\euffam=\seveneuf
\scriptscriptfont\euffam=\fiveeuf


\font\tenmsx=msam10
\font\sevenmsx=msam7
\font\fivemsx=msam5
\font\tenmsy=msbm10
\font\sevenmsy=msbm7
\font\fivemsy=msbm5
\newfam\msxfam
\newfam\msyfam
\textfont\msxfam=\tenmsx  \scriptfont\msxfam=\sevenmsx
  \scriptscriptfont\msxfam=\fivemsx
\textfont\msyfam=\tenmsy  \scriptfont\msyfam=\sevenmsy
  \scriptscriptfont\msyfam=\fivemsy
\edef\msx{\hexnumber\msxfam}

\mathchardef\upharpoonright="0\msx16
\let\restriction=\upharpoonright
\def\Bbb#1{\tenmsy\fam\msyfam#1}

\def\restrict{{\restriction}}

\def\Smallskip{\vskip1.4truecm}
\def\Bigskip{\vskip2.2truecm}

\def\qed{{\vcenter{\hrule height.4pt \hbox{\vrule width.4pt height5pt
 \kern5pt \vrule width.4pt} \hrule height.4pt}}}
\def\notin{{\in}\kern-5.5pt / \kern1pt}

\def\DD{{\Bbb D}}

\def\KK{{\Bbb K}}
\def\LL{{\Bbb L}}

\def\MM{{\Bbb M}}

\def\NN{{\Bbb N}}

\def\PP{{\Bbb P}}
\def\QQ{{\Bbb Q}}

\def\ZZ{{\Bbb Z}}

\font\small=cmr8 scaled\magstep0

\font\capit=cmcsc10 scaled\magstep0
\font\capitg=cmcsc10 scaled\magstep1

\font\dunhg=cmr10 scaled\magstep1
\font\dunhgg=cmr10 scaled\magstep2

\font\sanse=cmss10 scaled\magstep0

\font\ital=cmu10 scaled\magstep0
\overfullrule=0pt
\openup1.5\jot

\centerline{\dunhgg Evasion and prediction ---}
\Smallskip
\centerline{\dunhgg --- the Specker phenomenon and Gross spaces}
\Bigskip
\centerline{\capitg 
J\"org Brendle\footnote{$^*$}{{\small
The author would like to thank the MINERVA-foundation
for supporting him}}}
\Smallskip
{\baselineskip=0pt {\small
\noindent 
Department of Mathematics,
Bar--Ilan University,
52900 Ramat--Gan, Israel 
\smallskip
and
\smallskip
\noindent Mathematisches Institut der Universit\"at T\"ubingen,
Auf der Morgenstelle 10, 72076 T\"ubingen, Germany}}
\Smallskip
\centerline{\capit Abstract}
\medskip
\noindent We study the set--theoretic combinatorics underlying
the following two algebraic phenomena. \par
\item{(1)} A subgroup $G \leq \ZZ^\omega$ {\it exhibits
the Specker phenomenon} iff every homomorphism $G\to\ZZ$
maps almost all unit vectors to 0. Let {\bf se}
be the size of the smallest $G\leq\ZZ^\omega$ exhibiting
the Specker phenomenon. \par
\item{(2)} Given an uncountably dimensional vector space
$E$ equipped with a symmetric bilinear form $\Phi$ over
an at most countable field $\KK$, $(E,\Phi)$
is {\it strongly Gross} iff for all countably--dimensional $U \leq E$,
we have $dim(U^\perp) \leq \omega$.\par
\noindent Blass showed that the Specker phenomenon is closely
related
to a combinatorial phenomenon he called {\it evading and
predicting}.
We prove several additional results (both theorems of $ZFC$
and independence proofs) about evading and predicting as
well as {\bf se}, and relate a Luzin--style property associated
with evading to the existence of strong Gross spaces.
\Bigskip
{\dunhg Introduction}
\Smallskip
{\ital History and motivation.} The goal of this work is
the investigation of the set--theoretic combinatorics underlying two
algebraic phenomena which do not seem related at first
glance.
\par
One of them, coming from abelian group theory,
has been studied
in recent work of Blass [Bl 2]. 
We say that $G\leq(\ZZ^\omega,+)$ {\it exhibits the Specker
phenomenon} iff, given a group homomorphism $h: G \to
\ZZ$, for all but finitely many unit
vectors $e_n$ (i.e.
$e_n \in \ZZ^\omega$ is defined by $e_n (n) = 1$ and $e_n (i) = 0$
for $i \neq n$), one has $h (e_n) =0$. 
We let
\smallskip
\centerline{${\bf se} := min \{ \vert G\vert ;\;
\oplus_n \langle e_n \rangle \leq G \leq \ZZ^\omega \;\land\;
\forall$ homomorphisms $h:G \to \ZZ$, $h(e_n) = 0$
for almost all $n\}$,}
\smallskip
\noindent the {\it Specker--Eda number}.
Specker [S, Satz III] proved that $\ZZ^\omega$ exhibits
the Specker phenomenon, and Eda [E] showed that ${\bf
se}=2^\omega$ under Martin's Axiom $MA$ as well as the consistency of 
${\bf se} < 2^\omega$. Thus the question --- investigated
in [Bl 2] --- comes up whether any of the combinatorial cardinal
invariants of the continuum equals {\bf se}.
\par
The other phenomenon comes from quadratic form theory,
and has been intensively discussed in recent work of
Spinas, Baumgartner, and Shelah (see [BSp], [Sp 1],
[ShSp] and the survey [Sp 2]). Given a finite or countable
field $\KK$, an uncountably--dimensional $\KK$--vector
space $E$ equipped with a symmetric bilinear form $\Phi
: E^2 \to \KK$ is called {\it $($strongly$)$ Gross}
iff for all countably--dimensional subspaces $U \leq
E$, one has $dim U^\perp < dim E$ ($dim U^\perp \leq
\omega$) --- here, $U^\perp := \{ e \in E ; \; \forall u \in
U \; (\Phi (e,u) = 0 ) \}$ denotes the orthogonal complement
of $U$ in $E$. Many of the results in the work mentioned
above say roughly that, for certain cardinal invariants
{\bf ci} of the continuum, {\it if ${\bf ci} =
\omega_1$, then there is a strong Gross space over $\KK$},
and that for certain (different!) cardinals ${\bf ci'}$,
{\it if ${\bf ci'} > \omega_1$, then there are no strong Gross spaces
over $\KK$}. Thus the question naturally arises whether 
there is a cardinal (defined in terms of set--theoretic
combinatorics) such that its being $\omega_1$ is equivalent to
the existence of strong Gross spaces.
\par
Both phenomena turn out to be related to the following
combinatorial concept which has been introduced as well
in Blass' work [Bl 2, section 4]. Given a finite or
countable set $S$, an {\it $S$--valued predictor}
(or simply: {\it $S$--predictor}) is a pair $\pi
=(D_\pi , (\pi_n ; \; n \in D_\pi ))$ where $D_\pi
\subseteq \omega$ is infinite and for each $n \in D_\pi$,
$\pi_n$ is a function from $S^n $ to $S$. We say $\pi$
{\it predicts} $f \in S^\omega$ iff for all but finitely many
$n \in D_\pi$, we have $f(n) = \pi_n (f \restrict n)$; otherwise
$f$ {\it evades} $\pi$. We set
\smallskip
\centerline{${\bf e}:= min\{ \vert F \vert ;\; F \subseteq
\omega^\omega \;\land\; \forall \;\omega$--predictors $\pi
\; \exists f \in F \; (f$ evades $\pi ) \}$,}
\smallskip
\noindent the {\it evasion number}. A $\ZZ$--predictor is
called {\it linear} iff for all $n \in D_\pi$, the
function $\pi_n : \ZZ^n \to \QQ$ is linear with
rational coefficients. We let
\smallskip
\centerline{${\bf e}_\ell := min \{\vert F \vert ;\; F \subseteq
\ZZ^\omega \;\land\;\forall$ linear $\ZZ$--predictors $\pi \;
\exists f \in F \; (f$ evades $\pi ) \}$,}
\smallskip
\noindent the {\it linear evasion number}. Clearly
${\bf e}_\ell \leq {\bf e}$.
\bigskip
{\ital Results.}
To be able to explain our main results, we shall need
the definition of some of the classical cardinals
associated with the continuum --- for more on such
cardinals as well as the forcing notions used in our proofs
we refer the reader to section 1. Given a $\sigma$--ideal
${\cal I}$ on the real line, the {\it additivity of ${\cal I}$},
$add ({\cal I})$, is the smallest size of a family of members
of ${\cal I}$ the union of which is not in ${\cal I}$.
${\cal L}$ and ${\cal M}$ denote the ideals of
{\it Lebesgue measure zero} and {\it meager} sets,
respectively. For $A, B \subseteq \omega$,
$A \subseteq^* B$ ({\it $A$ is almost contained in
$B$}) means that $A \setminus B$ is finite;
and for $f , g \in\omega^\omega$, we say $f \leq^* g$ 
({\it $g$ eventually dominates $f$}) iff
$\forall^\infty n \; (f(n) \leq g(n))$ (where $\forall^\infty
n$ stands for {\it for all but finitely many $n$};
similarly $\exists^\infty n$ denotes {\it there are infinitely
many $n$}).
Using this notation, the {\it pseudointersection number}
{\bf p} is the smallest cardinality of a family $F$ of
subsets of $\omega$ with the {\it strong finite intersection
property} (i.e. given finitely many $A_i \in F$, $i<n$,
we have $\vert \bigcap_{i<n} A_i \vert = \omega$)
so that $\neg\exists B \in [\omega]^\omega \; \forall A
\in F \; (B \subseteq^* A)$ (we say: $F$ does not
have a {\it pseudointersection}). Next, the {\it
unbounding number} {\bf b} is the smallest size of a
family $F \subseteq \omega^\omega$ so that $\forall g
\in \omega^\omega \;\exists f \in F\; \exists^\infty
n \; (g(n) \leq f(n))$.
Finally, the {\it splitting number} {\bf s} is the smallest
cardinality of a family $F$ of subsets of $\omega$ so
that $\forall B \in [\omega]^\omega \;\exists A \in F \;
(\vert A\cap B \vert = \vert B\setminus A\vert =\omega )$
(we say: $A$ {\it splits} $B$).
It is well--known ([Fr], [vD, section 3]) 
that $\omega_1 \leq add({\cal L}),
{\bf p} \leq add({\cal M}) \leq {\bf b} \leq 2^\omega$
and ${\bf p} \leq {\bf s }\leq 2^\omega$.
\par
Blass [Bl 2, Theorem 12, Corollaries 11 and 14] proved
that $add ({\cal L}) \leq {\bf e}_\ell \leq
add ({\cal M})$, as well as ${\bf p} \leq {\bf e}_\ell$.
It is well--known that using standard techniques
(see 3.1. for details), the consistency of
${\bf p}, add({\cal L}) < {\bf e}_\ell$ with $ZFC$ can
be shown. We complete this cycle of independence results
(in 2.1. and 2.2.) by proving the
consistency of ${\bf e}_\ell < add({\cal M})$; more
explicitly:
\smallskip
{\capit Theorem A.} {\it It is consistent that ${\bf e}
=\omega_1 < {\bf b} = add({\cal M}) = 2^\omega = \kappa$
for any regular uncountable $\kappa$.}
\smallskip
\noindent Another result of Blass' concerns the relationship
between the Specker--Eda number and the concept of evasion,
namely ${\bf e}_\ell \leq {\bf se}$ [Bl 2, Corollary 8
and Theorem 10]. We shall see (in 2.4.) that an upper
bound to {\bf se} can be given in terms of evasion
as well (the cardinal ${\bf e'}$, introduced in 2.3.),
and derive from this:
\smallskip
{\capit Theorem B.} {\it ${\bf se} \leq unif({\cal L})$,
the size of the smallest non--measurable set of reals;
in particular ${\bf se} < {\bf b}$ is consistent.}
\smallskip
\noindent The interest in the latter consistency
stems from Blass' ${\bf se} \leq {\bf b}$ [Bl 2,
Theorem 2].
\par
In the third section we look at the phenomenon of
{\sanse evading and predicting} in general and in particular
at the relation between various forms of evasion
numbers (and some other cardinals as well).
Namely, we consider the spaces $n^\omega$ ($n \geq 2$),
$n$--valued predictors, and the corresponding evasion
numbers ${\bf e}_n$ --- or, more generally, 
compact spaces of the form $\prod_{n\in\omega}
f(n) = \{ g \in \omega^\omega ;\; \forall n\;
(g(n) < f(n)) \}$ for $f\in\omega^\omega$, and the
corresponding predictors and evasion numbers ${\bf e}_f$
(see 3.1. for exact definitions). We let
${\bf e}_{ubd} := min \{ {\bf e}_f ;\;
f\in\omega^\omega \}$ and
${\bf e}_{fin} := min \{ {\bf e}_n ;\;
n \in \omega \}$.
\noindent We shall show in 3.2. and 3.3.:
\smallskip
{\capit Theorem C.} {\it (a) ${\bf e} \geq min \{
{\bf b} , {\bf e}_{ubd} \}$ and ${\bf s} \leq {\bf e}_{fin}
= {\bf e}_n$ for all $n$. \par
\noindent (b) Both ${\bf e} < {\bf e}_{ubd}$ and ${\bf e}_{ubd}
< {\bf e}_{fin}$ are consistent.}
\smallskip
The forth section deals with Luzin--style properties
related to evading: given an arbitrary finite or
countable field $\KK$, an uncountably--dimensional 
subspace $(G,+) \leq \KK^\omega$
is called a {\it Luzin group} iff for all linear 
$\KK$--predictors $\pi$ all but countably many elements of
$G$ evade $\pi$; $G$ is {\it generalized Luzin} of size
$\kappa$ iff any linear $\KK$--predictor predicts
less than $\kappa$ many elements of $G$. We shall prove
(in 4.3.--4.6.):
\smallskip
{\capit Theorem D.} (dichotomy theorem) {\it 
(a) It is consistent that there are no generalized
Luzin groups $G \leq \KK^\omega$, where $\KK$ is any finite
field. \par
\noindent (b) 
For any countable field $\KK$, there is a generalized Luzin
group $G \leq \KK^\omega$ of size ${\bf b}$.}
\smallskip
{\capit Theorem E.} (equivalence theorem) {\it
For any finite or countable field $\KK$, the following are
equivalent: \par
\item{(a)} there exists a strong Gross space $(E, \Phi)$
over $\KK$; \par
\item{(b)} there is a Luzin group $G \leq \KK^\omega$.}
\smallskip
\noindent Using these results one gets alternative proofs
of the theorems of Baumgartner, Shelah and Spinas ([BSp],
[ShSp]) as well as one new result: $add({\cal L})
> \omega_1$ implies the non--existence of strong Gross
spaces (4.7.).
\par
We close our considerations with a list of questions in
section 5.
\bigskip
{\ital Notation.} We use standard set--theoretic notation
and refer the reader to [Ku], [Je 1], [Je 2] and
[Bau] 
for set theory in general
and forcing in particular. 
\par
Given a finite sequence $s$ (i.e. 
$s \in \omega^{<\omega}$), we let $lh(s) 
:= dom (s) $ denote the length of $s$; for $\ell \in lh(s)$,
$s \restrict \ell$ is the restriction of $s$ to $\ell$. 
$\hat{\;}$ is used for concatenation of sequences; and
$\langle \rangle$ is the empty sequence. 
--- Given a finite set $A \subseteq \kappa$ and $i < \vert A \vert$,
$A(i)$ denotes the $i$--th element of $A$ under the inherited
ordering.
--- 
Given a p.o. $\PP $, we
shall denote $\PP$-names by symbols like
$\breve f$, $\breve \pi$, $\breve D$, ...
\bigskip
{\ital Acknowledgment.} I would like to thank Myriam for drawing
the diagram on the computer.
\Bigskip
{\dunhg $\S$ 1. Cardinals and forcing notions}
\Smallskip
{\sanse 1.1. Cardinals.} In addition to the cardinals we have
seen already, we define, for a given $\sigma$--ideal ${\cal I}$
on the reals,
\smallskip
\item{} $cov({\cal I}) := min \{ \vert F\vert ;\; F \subseteq {\cal I}
\;\land\; \cup F = \omega^\omega \}$, the {\it covering
number} of ${\cal I}$,
\par
\item{} $unif({\cal I}) := min \{ \vert F\vert ;\; F \in
P(\omega^\omega) \setminus {\cal I} \}$, the {\it
uniformity} of ${\cal I}$, and
\par
\item{} $cof({\cal I}) := min \{ \vert F\vert ;\; F \subseteq 
{\cal I} \;\land\; \forall A \in {\cal I} \;\exists
B \in F \; (A\subseteq B) \}$, the {\it cofinality} of
${\cal I}$.
\smallskip
\noindent Furthermore we set
\smallskip
\item{} ${\bf d}:=min \{ \vert F\vert ;\; F \subseteq 
\omega^\omega \;\land\; \forall g\in\omega^\omega\;\exists f \in F
\; (g \leq^* f) \}$, the {\it dominating number}, and
\par
\item{} ${\bf r}:= min \{\vert F\vert ;\; F\subseteq [\omega]^\omega
\;\land\; \forall A\in [\omega]^\omega \;\exists B \in F\;
(\vert A\cap B\vert < \omega$ or $\vert B \setminus A\vert
< \omega )\}$, the {\it reaping number}.
\smallskip
\noindent Most of these invariants come in pairs, i.e. one
of them
can be defined from the other essentially by taking
negation and modifying the range of quantifiers. Compare, e.g.,
$cov({\cal I}) = min \{ \vert F \vert ;\; F \subseteq {\cal I}
\;\land\; \forall x \in \omega^\omega \; \exists A \in F
\; (x \in A) \}$ with $unif({\cal I}) = min \{
\vert F \vert ;\; F \subseteq \omega^\omega \;\land\;
\forall A \in {\cal I}\;\exists x \in F \;
(x \not\in A) \}$, or ${\bf b} = min \{ \vert F\vert ;\;
F \subseteq \omega^\omega \;\land\; \forall f \in \omega^\omega
\;\exists f \in F \; (g \not\leq^* f) \}$ with {\bf d}.
Other pairs are $(add({\cal I}), cof({\cal I}))$
and
$({\bf s}, {\bf r})$. One effect of this duality is 
that $ZFC$--proofs of inequalities between cardinals
dualize (e.g., ${\bf b } \geq add({\cal M})$ is
proved the same way as ${\bf d} \leq cof({\cal M})$);
another, that consistency proofs involving finite support iterations
dualize as well (see [Bl 1, in particular
section 5], [BaJS, section 1] or [Br]
for duality). The inequalities between these
cardinals (as well as some others which are crucial
for our investigations) which are provable in
$ZFC$ are displayed in the diagram in subsection 3.5.
\bigskip
{\sanse 1.2. Forcing notions.}
{\ital Hechler forcing.}
The {\it Hechler p.o.} $\DD$ is defined as follows:
\smallskip
\centerline{$(s,f) \in \DD \Longleftrightarrow s \in \omega^{< \omega}
\; \land \; f \in \omega^\omega \; \land \; s \subseteq f \; 
\land \; f$ strictly increasing}
\smallskip
\centerline{$(s,f)  \leq (t,g) \Longleftrightarrow s \supseteq t \;
\land \;
\forall n \in \omega \; (f(n) \geq g(n))$}
\smallskip
\noindent Following
Baumgartner and Dordal ([BD, $\S$ 2];
see also [BrJS, $\S$ 1]),
given $t \in \omega^{< \omega}$ strictly increasing 
and $A \subseteq \omega^{<
\omega}$, we define by induction when the {\it rank}
$rk (t, A)$ is $\alpha$. \par
\item{(a)} $rk(t,A) =0$ iff $t \in A$. \par
\item{(b)} $rk (t,A) = \alpha$ iff  for no $\beta < \alpha$
we have $rk (t,A) = \beta$, but there are $m \in \omega$ and
$\langle t_k ; \; k \in \omega \rangle$ such that $\forall
k \in \omega$: $t \subseteq t_k$, $t_k \in \omega^m$,
$t_k (lh(t)) \geq k$, and $rk(t_k , A) < \alpha$.  \par
\noindent Clearly, the rank is either $< \omega_1$ or undefined
(in which case we say $rk = \infty$). 
The following result is the main tool in the proof
of Theorem A (see 2.1.).
\smallskip
{\capit Lemma} (Baumgartner--Dordal [BD, $\S$ 2],
see also [BrJS, 1.2.])
{\it Let $I \subseteq \DD$ be dense. Set $A: = \{ t ; \;
\exists f \in \omega^\omega$ such that $(t,f) \in I \}$.
Then $rk(t^*,A) < \omega_1$ for any $t^* \in \omega^{< \omega}$.}
$\qed$
\smallskip
{\ital Mathias forcing.} The {\it Mathias p.o.} $\MM$
is defined as follows [Je 2, part one, section 3]:
\smallskip
\centerline{$(s,S) \in \MM \Longleftrightarrow s \in
\omega^{<\omega}\;\land\;S\in[\omega]^\omega\;\land\;
s$ strictly increasing $\land\; max(ran(s)) < min (S)$}
\smallskip
\centerline{$(s,S) \leq (t,T) \Longleftrightarrow s\supseteq t
\;\land\; S \subseteq T \;\land\; \forall i \in (lh(s)
\setminus lh(t))\;(s(i) \in T)$}
\smallskip
{\ital Laver forcing.} The {\it Laver p.o.} $\LL$ is defined
as follows [Je 2, part one, section 3]:
\smallskip
\centerline{$T\in\LL \Longleftrightarrow T\subseteq \omega^{<\omega}
$ is a tree $\land\;\exists\rho\in T\; \forall\sigma\in T \;(
\sigma \subseteq \rho \;\lor\;[\rho\subseteq\sigma \;\land\;
\exists^\infty n\;(\sigma\hat{\;}\langle n \rangle \in T)])$}
\smallskip
\centerline{$T\leq S \Longleftrightarrow T \subseteq S$}
\smallskip
\noindent The $\rho$ required to exist in the above definition
is usually called the {\it stem} of $T$, $stem(T)$. Furthermore,
for $\rho\in T$ we let $succ_T (\rho) := \{ n \in \omega ; \;
\rho\hat{\;}\langle n \rangle \in T \}$, the set of 
{\it successors}
of $\rho$ in $T$.
\smallskip
{\ital Laver property.} Both Laver and Mathias forcing 
as well as their countable support iterations have the
following property of p.o.'s $\PP$, sometimes referred
to as {\it Laver property}. Given $p\in\PP$, a function
$f\in\omega^\omega$ and a $\PP$--name $\breve g$ for an element
of $\prod_n f(n)$, there is $\phi \in \prod_n [f(n)]^n$
and $q\leq p$ so that
\smallskip
\centerline{$q\forces_\PP"\forall n \; (\breve g (n)
\in \phi (n))"$.}
\smallskip
\noindent See, e.g., [Bau, section 9] for details. One consequence
of this is that $\PP$ adds neither random nor Cohen reals.
\Bigskip
{\dunhg $\S$ 2. Evasion and the Specker phenomenon}
\Smallskip
{\sanse 2.1.} To prove the consistency of
${\bf e} < add({\cal M})$ we shall iterate Hechler forcing
$\kappa$ times with finite support over a model $V$
satisfying $CH$. To make the
argument that ${\bf e}$ is still $\omega_1$ at the end go through smoothly,
we shall consider the following property of p.o.'s $\PP$:
\smallskip
\item{$(**)$} given $F \subseteq \omega^\omega \cap V$,
$F\in V$, a family of functions below the identity
(i.e. $\forall n \;\forall f \in F \; (f(n) \leq n))$,
such that for any countable family of predictors $\Pi$
there is $f\in F$ evading all $\pi \in \Pi$, and $\langle
\breve \pi_n ; \; n\in\omega\rangle$ a sequence of $\PP$--names
for predictors, we can find a sequence $\langle \pi_n ;\;
n \in \omega\rangle \in V$ of predictors such that
whenever $f\in F$ evades all $\pi_n$, then
\par
\centerline{$\forces_\PP "f$ evades all $\breve\pi_n "$.}
\smallskip
{\capit Theorem.} {\it $\DD$ satisfies $(**)$.}
\smallskip
{\it Proof.} 
We shall use the notion of rank for $\DD$ as explained in
1.2. \par
Let $F$ be a family satisfying the requirements in
the definition of $(**)$, and let $\langle \breve\pi_n ;\; n\in\omega
\rangle$ be a sequence of $\DD$--names for predictors.
Associated with the name $\breve \pi_n$ we have the name $\breve D_n$
and the sequence of names $\langle \breve \pi_n^m ; \; m \in \omega
\rangle$ such that 
\smallskip
\centerline{$\forces_\DD "\breve \pi_n$ predicts on the set $\breve D_n$;
$\breve \pi^m_n$ is the predicting function on the $m$--th element
of $\breve D_n$."}
\smallskip
\noindent Fix $n, m$. Let $I^m_n: = \{ (t,f) ; \; (t,f)$ decides the
$m$--th element of $\breve D_n$ (and all preceding ones), say:
$(t,f) \forces_\DD "k$ is the $m$--th element of $\breve D_n$";
and $(t,f)$ decides $\breve \pi^m_n$ (and all $\breve \pi^i_n$, $i<m$)
on all sequences of length $k$ below the identity $\}$. All $I^m_n$
are dense and $I^{m+1}_n \subseteq I^m_n$. Thus, if $A^m_n :=
\{ t; \; \exists f \in \omega^\omega \; ((t,f) \in I^m_n) \}$, we
can apply Lemma 1.2. to $A^m_n$;
i.e. $rk(t,A^m_n) < \infty$ for all $t \in \omega^{<\omega}$ strictly
increasing.
\par
Now we define by recursion on rank:
\par
\item{---} when $t$ is $(n,m)$--{\it happy}, when it is $(n,m)$--{\it
sad}, and when it is {\it minimal $(n,m)$--sad};
\par
\item{---} $k(t,n,m) \in \omega$ and $\pi(t,n,m)$, where
$dom(\pi(t,n,m)) = \{ \sigma \in \omega^{< \omega} ; \;
lh(\sigma) = k(t,n,m)$ and $\forall i \in k(t,n,m) \; (\sigma
(i) \leq i) \}$ and $ran(\pi (t,n,m)) \subseteq k(t,n,m) + 1$,
for $t$'s which are $(n,m)$--happy;
\par
\item{---} sets $\tilde D (t,n,m)$ and predictors $\tilde\pi
(t,n,m)$ for $t$ which are minimal $(n,m)$--sad. \par
\noindent $\rightarrow$ $rk(t,A^m_n) = 0$. \par
\noindent Then we say $t$ is $(n,m)$--happy.
We choose $f$ such that $(t,f) \in I^m_n$. Let $k(t,n,m)$ be such
that 
\smallskip
\centerline{$(t,f) \forces_\DD "k(t,n,m)$ is the $m$--th element of
$\breve D_n$".}
\smallskip
\noindent Let $\pi(t,n,m)$ be such that, for $\sigma$ with $\forall
i \in k(t,n,m) \; (\sigma(i) \leq i)$,
$$\pi(t,n,m) (\sigma) = j \longleftrightarrow \cases{(t,f)
\forces_\DD "\breve \pi^m_n (\sigma) = j" &
and $j\leq k(t,n,m)$ or \cr
(t,f) \forces_\DD "\breve \pi^m_n (\sigma) >
k(t,n,m) " & and $j=0$ \cr }$$
This makes sense by definition of the set $I^m_n$.
\par
\noindent $\rightarrow$ $rk(t,A^m_n) = \alpha$. \par
\noindent Then we have $\langle t_i ; i \in \omega \rangle$ and
$\ell \in \omega$ such that for all $i$: $t \subseteq t_i$,
$t_i \in \omega^\ell$, $t_i (lh (t)) \geq i$ and $rk(t_i
, A^m_n ) < \alpha$. If $\exists^\infty i$ such that
$t_i$ is $(n,m)$--sad, then $t$ is $(n,m)$--sad, but
not minimal. Now suppose that $\forall^\infty i$, $t_i$ is
$(n,m)$--happy. \par
\noindent If $\{ k (t_i , n,m) ; \; i$ such that $t_i$ is
$(n,m)$--happy $\}$ is infinite, $t$ is minimal $(n,m)$--sad.
In this case we can without loss assume that $i<j$ implies
$k(t_i ,n,m) < k(t_j ,n,m)$. Let $\tilde D (t,n,m) =
\{ k (t_i ,n,m); \; i \in \omega \}$ and define
a predictor $\tilde\pi (t,n,m)$ as follows:
\par
\item{$\bullet$} $D_{\tilde\pi (t,n,m)} = \tilde D (t,
n,m)$; \par
\item{$\bullet$} for $k \in \tilde D (t,n,m)$, let
$i \in \omega$ be such that $k = k(t_i , n ,m)$ and
set $\tilde\pi_k (t,n,m) (\sigma) := \pi (t_i,n,m)
(\sigma)$ for $\sigma$ of length $k$ below the identity;
for other $\sigma$, $\tilde\pi_k (t,n,m)
(\sigma)$ can be defined arbitrarily.
\par
\noindent If $\{ k (t_i ,n,m) ; \; i$ such that $t_i$ is $(n,m)$--happy
$\}$ is finite, $t$ is still $(n,m)$--happy. In this case
we can without loss assume that the latter set contains
just one element, $k(t,n,m)$; and that $\forall i \in \omega$,
$\pi (t_i ,n,m)$ is the same function which we
call $\pi(t,n,m)$. --- This concludes the definition
of happiness and sadness. \par
Next, for each $n \in \omega$ and each $t$ such that
$\forall^\infty m$ $t$ is $(n,m)$--happy, we define
a predictor $\hat \pi (t,n)$ as follows:
\par
\item{$\bullet$} $D_{\hat \pi (t,n)} = \{ k(t,n,m) ;
\; m \in \omega \land t $ is $(n,m)$--happy $\}$
(note that this set must be infinite as $k(t,n,m) \geq m$);
\par
\item{$\bullet$} for $k \in D_{\hat\pi (t,n)}$, let $m \in
\omega$ be minimal such that $k=k(t,n,m)$ and set
$\hat\pi_k (t,n) (\sigma) := \pi(t,n,m) (\sigma)$
for $\sigma$ of length $k$ below the identity (for other
$\sigma$, $\hat\pi_k(t,n) (\sigma)$ can be defined
arbitrarily).
\par
\noindent Let $\Pi := \{ \tilde \pi (t,n,m) ; \; t$ minimal 
$(n,m)$--sad $\} \cup \{ \hat\pi (t,n) ; \; \forall^\infty
m$, $t$ is $(n,m)$--happy $\}$.
This is a countable set of predictors. 
Choose $f\in F$ evading all $\pi\in\Pi$.
\smallskip
{\it Claim. $\forces_\DD "f$ evades all $\breve\pi_n$, $n \in
\omega "$.}
\smallskip
{\it Proof of Claim.} By contradiction. Suppose there are
$(t,g) \in \DD$, $n \in \omega$ and $m_0 \in \omega$ such that
\smallskip
(+)\centerline{$(t,g) \forces_\DD " \forall m \geq m_0
\; (\breve \pi^m_n (f \restrict \breve k_m) = f (\breve k_m)
)$",}
\smallskip
\noindent where
\smallskip
\centerline{$\forces_\DD "\breve k_m$ is the $m$-th element
of $\breve D_n$".}
\smallskip
\noindent We consider two cases. \par
1. $\forall^\infty m \; (t$ is $(n,m)$--happy). \par
\noindent Then we look at the predictor $\hat \pi (t,n)$.
Let $m_1 \geq m_0$ be such that $\forall m \geq m_1$,
$t$ is $(n,m)$--happy. As $f$ evades $\hat\pi (t,n)$,
there is $m_2 \geq m_1$ and $k' \in D_{\hat \pi (t,n)}$
such that $m_2$ is minimal with $k' = k(t,n,m_2)$ and
$\hat\pi_{k'} (t,n) (f \restrict k') \neq f(k')$.
By construction, $\hat\pi_{k'} (t,n) (f \restrict k')
=\pi (t,n,m_2) (f\restrict k')$.
\par
{\it Subclaim 1. There is $(t',g') \leq (t,g)$ such
that  \smallskip
\centerline{$(t',g') \forces_\DD " k' = \breve k_{m_2}$ and
$\breve\pi^{m_2}_n (f \restrict k') = \pi (t,n,m_2) (f \restrict k')
"$,}
\smallskip
\noindent contradicting {\rm (+)}.}
\smallskip
{\it Proof of Subclaim 1.} This is an easy induction on
rank. If $rk(t,A^{m_2}_n ) = 0$, let $t'=t$ and $g' \geq g$
such that $(t,g') \in I^{m_2}_n$. If $rk(t, A^{m_2}_n) > 0$,
find $s$ such that $rk(s ,A^{m_2}_n ) < rk (t , A^{m_2}_n )$,
$s \supseteq t$, $\forall i \in dom (s) \; (s(i) \geq
g(i))$, $s$ is $(n,m_2)$--happy,
$k(s,n,m_2) = k(t,n,m_2) = k'$, and $\pi(s,n,m_2) =
\pi(t,n,m_2)$. $\qed$
\smallskip
2. $\exists^\infty m \; (t$ is $(n,m)$--sad).
\par
\noindent Choose $m_1 \geq m_0$ such that $t$ is $(n,m_1)$--sad.
Next choose $(t',g') \leq (t,g)$ such that $t'$ is minimal
$(n,m_1)$--sad (this is possible by construction). This time
we look at the predictor $\tilde\pi (t',n,m_1)$. Choose $i_0$
such that $\forall i \geq i_0 \; \forall j \in dom (t_i)
\; (t_i (j) \geq g' (j))$, where the sequence $\langle t_i;
\; i \in \omega \rangle$ is chosen for $t$ as in the definition
of minimal $(n,m_1)$--sadness. As $f$ evades $\tilde \pi
(t',n,m_1)$, there is $k' \in \tilde D (t', n,m_1)$,
$k' \geq k(t_{i_0},n,m_1)$, such that $\tilde\pi_{k'}
(t',n,m_1) (f \restrict k') \neq f(k')$. By construction,
$\tilde\pi_{k'} (t',n,m_1) (f\restrict k') = \pi (t_i , n,m_1)
(f \restrict k')$, where $i \geq i_0$
is such that $k' = k(t_i , n,m_1)$.
\smallskip
{\it Subclaim 2. There is $(t'', g'') \leq (t',g')$
such that
\smallskip
\centerline{$(t'',g'') \forces_\DD "k' = \breve k_{m_1} \land
\breve\pi^{m_1}_n (f \restrict k') = \pi (t_i ,n,m_1)
(f\restrict k')"$,}
\smallskip
\noindent contradicting $(+)$. }
\smallskip
{\it Proof of Subclaim 2.} Again an easy induction on rank. If
$rk(t_i , A^{m_1}_n) = 0$, let $t'' = t_i$ and $g'' \geq g'$
such that $(t_i , g'') \in I^{m_1}_n$. Then clearly 
$(t'',g'') \leq (t',g')$.
If $rk(t_i , A^{m_1}_n ) >0$, we proceed as in the proof
of subclaim 1. $\qed$
\smallskip
This concludes the proof of the claim and finishes the proof 
of the Theorem as well. $\qed$ $\qed$
\bigskip
{\sanse 2.2.} 
Now let $\DD_\alpha$ denote the iteration of
Hechler forcing of length $\alpha$. We claim that $\DD_\alpha$
still has property $(**)$. By 2.1. this
is a consequence of the following preservation result:
\smallskip
{\capit Lemma.} {\it Assume $\langle \PP_\beta , \breve\QQ_\beta ;
\;
\beta < \alpha \rangle$ is an $\alpha$--stage finite support
iteration of $ccc$ partial orders such that \smallskip
\centerline{$\forall \beta < \alpha \;\;\; \forces_{\PP_\beta}
"\breve \QQ_\beta$ satisfies $(**)"$.}
\smallskip
\noindent Then $\PP_\alpha$
satisfies $(**)$.}
\smallskip
{\it Proof.} By induction on $\alpha$. The successor step as
well as the case $cf(\alpha) > \omega$ are
trivial. So assume $cf(\alpha) = \omega$; without loss
$\alpha = \omega$.
\par
Let $\langle \breve\pi_n ; \; n \in \omega \rangle$ be a sequence
of $\PP_\omega$--names for predictors; for $n \in \omega$,
$\breve\pi_n = (\breve D_n ; \; (\breve\pi^k_n ; \;
k \in \omega ))$. For each $m \in \omega$ let $\langle
\breve\pi_{n,m} ; \; n \in \omega\rangle = \langle
(\breve D_{n,m} ; \; (\breve\pi^k_{n,m} ; \; k \in
\omega )); \; n \in \omega \rangle$ be a sequence of
$\PP_m$--names for predictors and $\langle \breve p^k_m ; \;
k \in \omega \rangle$ a sequence of $\PP_m$--names for
elements of $\breve \PP_{[m,\omega)}$ such that
\smallskip
\itemitem{$\forces_{\PP_m}$} $\breve p^{k+1}_m \leq
\breve p^k_m$ and $\breve p^k_m \forces_{\breve \PP_{[m,
\omega)}}"$the $k$--th elements of $\breve D_n$ and
$\breve D_{n,m}$ are equal, say $\ell$, and $\breve \pi^k_n
$ equals $\breve \pi^k_{n,m}$ on all sequences of length 
$\ell$ below the identity".
\smallskip
\noindent By induction hypothesis find $\langle \pi_{n,m}
;\; n,m \in \omega \rangle \in V$ a sequence of predictors such that
whenever $f \in F$ evades all $\pi_{n,m}$, then for
all $m$:
\smallskip
\centerline{$\forces_{\PP_m} "f$ evades all $\breve\pi_{n,m}$,
where $n \in \omega$".}
\smallskip
\noindent We claim that $\forces_{\PP_\omega} "f$
evades all $\breve\pi_n$".
\par
For suppose there are a condition $p \in \PP_\omega$,
$ n \in \omega$ and $k_0 \in \omega$ such that
\smallskip
(+)\centerline{$p \forces_{\PP_\omega} "\forall k \geq k_0
\; (\breve \pi^k_n (f\restrict \breve \ell^k) = f (\breve \ell^k)
)$",}
\smallskip
\noindent where $\forces_{\PP_\omega} "\breve\ell^k$
is the $k$-th element of $\breve D_n$". Let $m = supp(p)$.
By induction hypothesis we know that 
\smallskip
\centerline{$\forces_{\PP_m} "f$ evades $\breve \pi_{n,m}$".}
\smallskip
\noindent Hence we can find $q \leq p$, $q \in \PP_m$, and
$k \geq k_0$ such that 
\smallskip
\centerline{$q \forces_{\PP_m} "f(\breve \ell^k_m) \neq
\breve \pi^k_{n,m} (f \restrict \breve \ell^k_m)"$,}
\smallskip
\noindent where $\forces_{\PP_m} "\breve \ell^k_m$ is the $k$--th
element of $\breve D_{n,m}$". Thus, by definition of the
name $\breve p^k_m$,
\smallskip
\centerline{$q \forces_{\PP_m} \breve p^k_m \forces_{\breve \PP_{
[m,\omega)}} "f(\ell^k_m) \neq \breve \pi^k_{n,m} (f 
\restrict \breve \ell^k_m) = \breve \pi^k_n (f \restrict
\ell^k)$",}
\noindent contradicting (+). $\qed$
\smallskip
Applying $(**)$ to $\DD_\kappa$ we get that $F
=(\omega^\omega)^V$ is a family of functions of size $\omega_1$
such that for every predictor $\pi \in V^{\DD_\kappa}$,
there is $f\in F$ evading $\pi$. Thus $V^{\DD_\kappa}
\models {\bf e}= \omega_1$. It is well--known that $V^{\DD_\kappa}
\models add({\cal M}) = \kappa$.
This ends the proof of Theorem A. Using the methods of [Br]
one can in fact show the consistency of ${\bf e} = \kappa$
and ${\bf b} = add({\cal M}) =\lambda$ for any regular
uncountable $\kappa < \lambda$.
\bigskip
{\sanse 2.3.} We consider the following more general notion
of predicting: we are given two sets $D_\pi = \{ k_n ;
\; n\in\omega\}\subseteq\omega$ and $E_\pi =\{ \ell_n ;\; n\in
\omega\}\subseteq\omega$
such that $k_n \leq \ell_n < k_{n+1}$ for all $n \in \omega$;
we also have for each $n$ a function $\pi_n :\omega^{\ell_n
\setminus\{k_n\}} \to \omega$; we say the {\it predictor}
$\pi = (D_\pi , E_\pi , (\pi_n ; \; n\in\omega))$ {\it
predicts} $f\in\omega^\omega$ iff $\forall^\infty n
\;(\pi_n (f\restrict (\ell_n\setminus\{ k_n\})) =
f(k_n))$.
We let ${\bf e'}$ be the smallest size of a set of
functions $F$ from $\omega$ to $\omega$ such that given
a countable set of such predictors $\Pi$, there is
$f\in F$ evading all $\pi \in \Pi$.
Clearly, ${\bf e'} \geq {\bf e}$. Also, a set predicted
by countably many predictors is necessarily contained in the
union of countably many closed measure zero sets.
Thus ${\bf e'} \leq unif({\cal M}), unif({\cal L})$.
\bigskip
{\sanse 2.4.} We shall give an upper bound to ${\bf se}$
in terms of evading by showing:
\smallskip
{\capit Theorem.} {\it ${\bf se} \leq {\bf e'}$.}
\smallskip
\noindent Note that, by the remarks in 2.3., this finishes
the proof of Theorem B: to get the consistency of
${\bf se} < {\bf b}$ simply add $\omega_1$ random reals
over a model for $MA$; then ${\bf se} \leq {\bf e'}
\leq unif({\cal L}) = \omega_1$, whereas ${\bf b}
=2^\omega$.
\smallskip
{\it Proof of Theorem.} Let $\{ p_n ;\; n \in \omega \}$
be an enumeration of all primes. Let ${\cal F} =
\{ f_\alpha ; \;\alpha < {\bf e'} \} \subseteq \omega^\omega$
be a family of functions evading all families of countably
many predictors (in the sense of 2.3., of course).
For $\alpha < {\bf e'}$, we define $x_\alpha \in \omega^\omega$
as follows:
\smallskip
\centerline{$x_\alpha (0) = 1$}\par
\centerline{$x_\alpha (1) = p_0^{f_\alpha (0)}$}\par
\centerline{$x_\alpha (2) = p_0^{f_\alpha (0) + 1} \cdot
p_1^{f_\alpha (1)}$}\par
\centerline{$...$}\par
\centerline{$x_\alpha (n) = \prod_{i<n} p_i^{f_\alpha(i)
+n-i-1}$}
\smallskip
\noindent We let $G \leq \ZZ^\omega$ be the pure subgroup
of $\ZZ^\omega$ generated by the $x_\alpha$ and the unit
vectors $e_n$. Clearly $\vert G\vert = {\bf e'}$.
We claim that $G$ exhibits the Specker phenomenon.
\par
For suppose not, and assume $h: G \to \ZZ$ is a homomorphism
such that there are infinitely many $n$ with $h(e_n) \neq 0$.
Fix $z\in\ZZ$. We shall introduce a predictor
$\pi=\pi_z$ of the required sort. $D_\pi =
\{ k_n ;\; n \in \omega\}$ and $E_\pi =\{\ell_n ;\; n\in\omega\}$
are such that
\par
\item{---} $h(e_{k_n}) \neq 0$, $ \vert h(e_{k_n})
\vert > 2 \cdot\sum_{i<k_n} \vert h(e_i) \vert$, and
$\forall \alpha < {\bf e'} \; \forall^\infty n 
\; ( \vert h(e_{k_n}) \vert > f_\alpha (k_n), x_\alpha (k_n))$
[Note that we can indeed choose the $k_n$ in such a
way by the argument of the proof of [Bl 2, Theorem 2]
(going over to $G' \geq G$ still of size ${\bf e'}$,
if necessary).];
\par
\item{---} $k_0$ is such that $\vert h(e_{k_0}) \vert
> \vert z \vert$ and $k_n < \ell_n < k_{n+1}$ are such that
$\ell_n - k_n > 2 \cdot\vert h(e_{k_n}) \vert^2$.
\par
\noindent To motivate the predictor, assume there
are $f_\alpha, f_\beta \in {\cal F}$
such that $h(x_\alpha) = h(x_\beta) = z$,
$f_\alpha \restrict \ell_n \setminus \{k_n\}
=f_\beta\restrict\ell_n\setminus\{k_n\}$, and
$f_\alpha(k_n) \neq f_\beta (k_n)$,
where $n$ is large enough (so that
$\vert h(e_{k_n}) \vert > f_\gamma (k_n),
x_\gamma (k_n)$, where $\gamma = \alpha$ or
$\beta$). We put
$$a:= h(x_\alpha \restrict [e_{k_n+1} , \infty))
= z - \sum_{i<k_n} h(e_i) \cdot x_\alpha (i)
- h(e_{k_n}) \cdot x_\alpha (k_n) \neq 0,$$
because the absolute value of the last term is
larger than the absolute values of the others together.
Also note that $a=h(x_\beta \restrict [e_{k_n+1},
\infty ))$. For $j$ with $k_n +1 \leq j \leq \ell_n$,
we let $b_j = \prod_{i<k_n \lor
k_n<i<j} p_i^{f_\alpha (i) +j-i-1}$;
and, letting $p:=p_{k_n}$,
we define by recursion on such $j$ the numbers
$\hat x^\gamma_j$ where $\gamma = \alpha$ or $\beta$:
\smallskip
\centerline{$\hat x^\gamma_{k_n +1} \cdot b_{k_n + 1}\cdot
p^{f_\gamma (k_n)} = a = h(e_{k_n+1})\cdot b_{k_n+1}
\cdot p^{f_\gamma (k_n)} + \hat x^\gamma_{k_n +2}
\cdot b_{k_n +2} \cdot p^{f_\gamma (k_n) + 1}$}
\medskip \centerline{$\hat x^\gamma_j \cdot b_j \cdot
p^{f_\gamma (k_n)+j-k_n-1} = h(e_j) \cdot b_j \cdot p^{f_\gamma 
(k_n) +j-k_n-1}
+\hat x^\gamma_{j+1} \cdot b_{j+1} \cdot p^{f_\gamma (k_n)+j-k_n}$}
\smallskip
\noindent for $k_n +1<j<\ell_n$.
Note that, by purity, all the $\hat x^\gamma_j$
must be integers. Furthermore $\hat x^\gamma_{k_n + 1}
\neq 0$ because $a \neq 0$. 
Suppose without loss that $f_\alpha
(k_n) < f_\beta (k_n)$. Let $m$ be maximal so that
$\hat x^\alpha_{k_n +1}$
is divisible by $p^m$;
then $\hat x^\beta_{k_n + 1}$ is divisible only by $p^{m'}$
where $m' = m - f_\beta (k_n) + f_\alpha (k_n)$.
In particular $\hat x_{k_n+1}^\alpha - \hat x_{k_n+1}^\beta
\neq 0$.
By assumption on $n$, we must have $m < \ell_n - k_n$
(in fact, $\vert a \vert < \ell_n - k_n$).
We can now divide the above equations for $\alpha$ and
$\beta$, respectively, by the appropriate power of $p$
and products of other primes,
subtract recursively one from the other
and get
\smallskip
\centerline{$\hat x^\alpha_{k_n +1} - \hat x^\beta_{k_n +1}
=(\hat x^\alpha_{k_n + 2} - \hat x^\beta_{k_n + 2})
\cdot\prod_{i<k_n} p_i\cdot p_{k_n+1}^{f_\alpha (k_n +1)} \cdot p$}
\medskip\centerline{$\hat x^\alpha_j - \hat x^\beta_j
=(\hat x^\alpha_{j+1} - \hat x^\beta_{j+1})
\cdot\prod_{i<k_n\lor k_n<i<j} p_i\cdot p_j^{f_\alpha (j)} \cdot p$}
\smallskip
\noindent Then $\hat x^\alpha_{k_n+1} - \hat x^\beta_{k_n+1}
$ is at most divisible by $p^{m'}$, hence $\hat x^\alpha_{k_n+2}
-\hat x^\beta_{k_n +2}$ is at most divisible by
$p^{m'-1}$ etc. Thus there must be a $j \leq \ell_n$
so that $\hat x^\alpha_j - \hat x^\beta_j$ is not an
integer
anymore, a contradiction.
\par
This shows that given $f_\alpha , f_\beta \in {\cal F}$
with $h(x_\alpha) = h(x_\beta) = z$ and $n$ large
enough, if $f_\alpha \restrict (\ell_n \setminus
\{k_n\}) = f_\beta \restrict (\ell_n \setminus \{k_n\})$,
then we must have $f_\alpha (k_n) = f_\beta (k_n)$;
so we simply let the predictor $\pi_z$ predict this uniquely
defined value. In the end we get a countable family
$\{ \pi_z; \; z \in\ZZ\}$ of predictors so that each
$f_\alpha$ is predicted by one $\pi_z$, a contradiction
to the choice of the family ${\cal F}$. $\qed$
\Bigskip
{\dunhg $\S$ 3. Towards a general theory of evasion
and prediction}
\Smallskip
{\sanse 3.1.} The proofs of the preceding section suggest that
we generalize the notions of evading and predicting defined in
the Introduction, and look at the corresponding cardinals.
\par
One way of doing this goes as follows.
Fix $f \in (\omega + 1 \setminus 2 )^\omega$, and let $X :=
\prod_n f(n)$; i.e. $X$ consists of the functions from the Baire
space which are below $f$ everywhere. An {\it $X$--predictor}
(or: {\it $f$--predictor}) is a pair $\pi = (D_\pi , 
(\pi_n ; \; n \in D_\pi ))$ such that for every $n \in D_\pi$,
$\pi_n : \prod_{k<n} f(k) \to f(n)$; $\pi$ {\it predicts}
$g \in X$ iff $\forall^\infty n \in D_\pi \; (\pi_n (g \restrict
n ) = g(n))$; otherwise $g$ {\it evades} $\pi$.
Let ${\bf e}_X$ (or: ${\bf e}_f$) be the corresponding {\it
evasion number}; i.e. the smallest size of a set of functions $F 
\subseteq X$ such that every $X$--predictor is evaded by
some $g \in F$. In case $f $ eventually equals $\omega$,
we get ${\bf e}_f = {\bf e}$; in case $X = n^\omega$
($n \geq 2$), we talk about $n$--predictors (instead
of $n^\omega$--predictors) and set ${\bf e}_n := {\bf e}_X$.
Finally we let ${\bf e}_{ubd} := min \{ {\bf e}_f ; \;
f \in \omega^\omega \}$, the {\it unbounded evasion number},
and ${\bf e}_{bdd} = min \{ {\bf e}_f ; \; f \in \omega^\omega$
is bounded $\} = min \{ {\bf e}_n ; \; n \in \omega \} =:
{\bf e}_{fin}$, the {\it bounded} or {\it finite evasion
number}. We trivially have ${\bf e} \leq {\bf e}_{ubd} \leq
{\bf e}_{fin}$, and we shall see in the next two subsections
that nothing else can be proved in $ZFC$ about the relationship 
between these three cardinals.
Furthermore, a set predicted by any predictor is easily
seen to be the union of countably many closed measure
zero sets; thus ${\bf e}_{fin} \leq unif({\cal E})
\leq unif({\cal M}), unif({\cal L})$,
where ${\cal E}$ is the $\sigma$--ideal generated by
the latter (this has been studied recently in [BS]).
\par
The notion of {\it linear predicting} can be generalized as well.
Let $\KK$ be a finite or countable field. A $\KK$--predictor
$\pi = (D_\pi , (\pi_n ;\; n \in D_\pi ))$ is called
{\it linear} iff for every $n \in D_\pi$,
$\pi_n : \KK^n \to \KK$ is a linear function. We let ${\bf e}_\KK$
be the smallest size of a set of functions $F \subseteq \KK^\omega$
such that every linear $\KK$-predictor is evaded by some
$f \in F$. As it is easier to evade just linear predictors,
we have ${\bf e}_\KK \leq {\bf e}_{\vert \KK \vert}$,
in particular ${\bf e}_\KK \leq {\bf e}$ for countable fields.
If $\QQ$ denotes the field of the rationals, ${\bf e}_\QQ =
{\bf e}_\ell$, the cardinal studied by Blass [Bl 2, section 4].
His results are easily seen to carry over to the other ${\bf e}_\KK$'s,
and we have, e.g., ${\bf e}_\KK \geq add({\cal L}) , {\bf p}$ for
any finite or countable field $\KK$ and ${\bf e}_\KK
\leq add({\cal M})$ for any countable field $\KK$ (cf
Introduction).
\par
Turning to consistency results, the following are known.
Iterating a Borel $\sigma$--centered forcing notion adding
a generic linear predictor (see 3.3. for similar forcing notions)
with finite support, we easily get $CON({\bf e} , {\bf e}_\KK
> add({\cal L}) , {\bf s}, {\bf p})$ [${\bf e}_\KK  > add({\cal L})$
holds because iterations of $\sigma$--centered forcing notions
do not add random reals (as far as I know this is due to
Miller and implicit in [Mi 1, section 5]); ${\bf e}_\KK > 
{\bf s}, {\bf p}$ holds
by easy definability of the forcing notions (see [JS 1])].
We finally note that in the proof of Theorem A (2.1. and 2.2.) we proved
in fact that ${\bf e}_{id} = \omega_1$ in the Hechler real
model; however, the choice of the identity function
was arbitrary, and it is easily
read off from the proof that all evasion numbers defined above
equal $\omega_1$ in the latter model.
\bigskip
{\sanse 3.2.} {\it Some $ZFC$--results.} 
We start the proof of Theorem C with a series
of Lemmata.
\smallskip
{\capit Lemma 1.} {\it ${\bf e} \geq min \{ {\bf b} , {\bf e}_{ubd}
\}$, and thus $min \{ {\bf b}, {\bf e} \} = min \{ {\bf b}
, {\bf e}_{ubd} \}$.}
\smallskip
{\it Proof.} Let $F \subseteq \omega^\omega$ be of size $< min
\{ {\bf b} , {\bf e}_{ubd} \}$. We have to find a predictor
predicting every function in $F$.
\par
By $\vert F \vert < {\bf b}$ find $x \in \omega^\omega$
such that $\forall f \in F \; \forall^\infty n \;
(f(n) < x(n))$. Given $f \in F$, we let for $n \in
\omega$
$$f' (n) := \cases{ f(n) & if $f(n) < x(n)$ \cr
0 & otherwise, \cr}$$
and let $F' := \{ f' ; \; f \in F \}$. By $\vert F' \vert
< {\bf e}_X$ (where
$X = \prod_n x(n)$) find an $X$--predictor $\pi = (D_\pi , (\pi_k ;
\; k \in D_\pi ))$ predicting every
function of $F'$. We extend $\pi$ to an $\omega$--predictor
$\pi^*$ as follows:
given $k \in D_\pi$ and $\sigma \in \omega^k \setminus
\prod_{n<k} x(n)$, we define for $n < k$
$$\sigma ' (n) := \cases{\sigma (n) & if $\sigma (n) < x(n)$ \cr
0 & otherwise. \cr}$$
Next we let $\pi^*_k (\sigma) := \pi_k (\sigma ')$. We claim that
$\pi^* = (D_{\pi} , (\pi^*_k ; \; k \in D_{\pi} ))$
predicts every function of $F$. 
\par
For, given $f \in F$, there is $n_0$ such that $\forall
n \geq n_0 \; (f' (n) = f(n))$; next there is $n_1 \geq n_0$
such that for all $k \geq n_1$, if $k \in D_\pi$,
then $\pi_k (f' \restrict k ) = f' (k)$.
Thus for all $k \geq n_1$, if $k \in D_\pi$,
then $f(k) = f' (k) = \pi_k (f' \restrict k) =
\pi^*_k (f \restrict k)$. $\qed$
\smallskip
{\capit Lemma 2.} {\it ${\bf e}_{fin} = {\bf e}_n$ for all 
$n \geq 2$.}
\smallskip
{\it Proof.} It suffices to show by induction on $n$
that ${\bf e}_n \geq {\bf e}_2$. To this end 
assume ${\bf e}_n < {\bf e}_2$ ($n$ minimal), and let
$F = \{ f_\alpha ; \; \alpha < {\bf e}_n \} \subseteq n^\omega$
be a family of functions evading every $n$--predictor.
\par
We define for $\alpha < {\bf e}_n$
$$g_\alpha (m) := \cases{0 \longleftrightarrow f_\alpha (m) \neq
1 &\cr
1 \longleftrightarrow f_\alpha (m) = 1. &\cr}$$
By assumption, there is a $2$--predictor $\pi =
(D_\pi , (\pi_{k_m} ; \; m \in \omega ))$
(where $D_\pi = \{ k_m ; \; m \in \omega \}$ is the
increasing enumeration)
predicting all $g_\alpha$. Next define for $\alpha < {\bf e}_n$
$$h_\alpha (m) := \cases{0 \longleftrightarrow f_\alpha
(k_m) = 0 {\rm \; or \;} 1 &\cr
f_\alpha (k_m) -1 \longleftrightarrow f_\alpha (k_m) \geq 2. &\cr
}$$
By assumption again, there is an $(n-1)$--predictor
$\pi ' = (D_{\pi '} , (\pi_m ' ; \; m \in D_{\pi '} ))$
predicting all $h_\alpha$. We can define an $n$--predictor
$\tilde \pi$ as follows: $D_{\tilde\pi} = \{ k_m ;\;
m \in D_{\pi '} \}$; given $\sigma \in n^{k_m}$ we
let 
$$\sigma ' (i) := \cases{ 0 \longleftrightarrow \sigma (k_i)
= 0 {\rm \; or \;} 1 &\cr
\sigma(k_i) - 1 \longleftrightarrow \sigma (k_i) \geq 2; &\cr}
$$
next, given $m \in D_{\pi '}$, we define 
$$\tilde\pi_{k_m} (\sigma) := \cases{
1 \longleftrightarrow \pi_{k_m} (\sigma) = 1 \;\land\;
\pi_m ' (\sigma ') = 0 &\cr
\ell \longleftrightarrow \pi_{k_m} (\sigma) = 0
\;\land\; \pi_m ' (\sigma ') = \ell - 1 & for
$\ell \geq 2$ \cr
0 \longleftrightarrow {\rm otherwise} &\cr}$$
\par
It is easy to check that $\tilde\pi$ predicts all
$f_\alpha$, thus giving a contradiction. $\qed$
\smallskip
{\capit Lemma 3.} {\it ${\bf e}_{fin} \geq {\bf s}$;
also ${\bf e}_\KK \geq {\bf s}$ for all finite fields
$\KK$.}
\smallskip
{\it Proof.} As ${\bf e}_{fin} \geq {\bf e}_\KK$, it suffices
to show the second inequality. To this end, let $\kappa < {\bf s}$
and let $\{ f_\alpha ; \; \alpha < \kappa \}$ be a family
of functions from $\omega$ to $\KK$. For each
$a \in \KK$ let $A_{a,\alpha} : = f_\alpha^{-1} (\{ a \})$.
As $\{ A_{a,\alpha } ; \; a \in \KK \;\land\;
\alpha < \kappa \}$ is not a splitting family, there is an
infinite $B \subseteq \omega$ which is not split, i.e.
$\forall \alpha < \kappa \; \exists a_\alpha \in \KK$ such that $
B \subseteq^* A_{a_\alpha,\alpha}$. Assume $B = \{ b_n ; \; n \in \omega
\}$ is the increasing enumeration of $B$. Let $D_\pi
:= \{ b_{2n + 1} ;\; n \in \omega \}$ and define
$\pi_n : \KK^{b_{2n+1}} \to \KK$ by $\pi_n (\sigma) =
\sigma (b_{2n})$. $\pi_n$ is trivially linear,
and $\pi = (D_\pi , (\pi_n ; \; n \in \omega))$ is
easily seen to predict all the $f_\alpha$. $\qed$
\smallskip
Let us note that these results together with Blass'
$min\{ {\bf e} , {\bf b} \} \leq add ({\cal M})$
[Bl 2, Theorem 13] already yield one consistency result
concerning evasion numbers, namely $CON ({\bf e}_{fin}
> {\bf e}_{ubd} )$. To see that this holds in the
Mathias real model, note that the latter satisfies
${\bf s} = \omega_2$ (this is straightforward from the
combinatorial properties of the Mathias generic real)
and hence ${\bf e}_{fin}  = \omega_2$, whereas 
$min \{ {\bf e}_{ubd} , {\bf b} \} =
min \{ {\bf e} , {\bf b} \} \leq add ({\cal M}) 
\leq cov({\cal M}) = \omega_1$
by the fact that iterated Mathias forcing doesn't
add Cohen reals (see 1.2.). Thus
${\bf b} = \omega_2$ in the Mathias real model gives
${\bf e}_{ubd} = \omega_1$. In fact, a canonical application
of the Laver property gives ${\bf e}_x = \omega_1$ for
any $x \in \omega^\omega$ converging
to infinity in this model.
\bigskip
{\sanse 3.3.} The goal of this subsection is to show:
\smallskip
{\capit Theorem.} {\it For any regular uncountable $\kappa$, it
is consistent that ${\bf e}_{ubd} = 2^\omega = \kappa$ and
${\bf e}={\bf b} = \omega_1$.}
\smallskip
\noindent This will complete the proof of Theorem C.
\par
Given $x \in \omega^\omega$, the p.o. $\PP_x$ 
for adding a predictor for
the space $X = \prod_n x(n)$ is defined as follows.
$$(d, \langle \pi_k ;\; k\in d \rangle , F) \in \PP_x 
\longleftrightarrow \cases{ d \subseteq \omega {\rm \; finite}
\;\land &\cr
\forall k \in d \; (dom (\pi_k) = \prod_{n<k}
x(n) \;\land\; ran(\pi_k) \subseteq x(k) ) \;\land &\cr
F \; {\rm finite}, \; \forall f \in F \; (dom (f) \leq
\omega \;\land\; f \in \prod_{n < dom(f)} x(n) ) &\cr}$$
We order $\PP_x$ by setting $(e, \langle \tilde\pi_k
;\; k \in e \rangle , G ) \leq (d , \langle \pi_k
;\; k \in d \rangle , F)$ iff
\par
\item{1)} $e \supseteq d$ and $max(d) < min (e \setminus d)$;
\par
\item{2)} $\tilde\pi_k = \pi_k$ for $k \in d$;
\par
\item{3)} $\forall f \in F \; \exists g \in G \;
(f \subseteq g)$ and $\forall k \in (e \setminus d )\;
\forall f \in F \; (k \in dom (f) \longrightarrow
\tilde\pi_k (f \restrict k) = f(k) )$.
\par
\noindent Clearly $\PP_x$ is a $\sigma$--centered
p.o. (more explicitly, $(d, \langle \pi_k
;\; k \in d \rangle , F)$ and $(d , \langle \pi_k ;\;
k \in d \rangle , G)$ are compatible with 
common extension $(d , \langle \pi_k ; \; k \in d \rangle
, F \cup G )$ for any choice of $F$ and $G$).
\par
Obviously, we can make ${\bf e}_{ubd} = \kappa$ by starting
with a model for $CH$, iterating p.o.'s of the form
$\PP_x$, where $x$ is a real in some intermediate stage of
the extension, $\kappa$ many times with finite support
(this is a standard enumeration argument as in the classical
consistency proof of $MA$).
\par
So we will be done if we can show that in the final model,
the ground model reals are still unbounded, and there is no
predictor predicting all ground model 
functions (in $\omega^\omega$); i.e. ${\bf b} =
{\bf
e} = \omega_1$ (by $CH$).
\par
To this end we use a modification of a notion and some techniques
of [BrJ, $\S$ 1]. Given a p.o. $\PP$, a function
$h : \PP \to \omega$ is a {\it height function} iff
$p \leq q$ implies $h(p) \geq h(q)$ for $p , q \in \PP$.
A pair $(\PP , h)$ is {\it soft} iff $\PP$ is a p.o.,
$h$ is a height function on $\PP$, and the following
three conditions are met:
\par
\item{(I)} if $\{ p_n ; \; n \in 
\omega \}$ is decreasing and $\exists m \in \omega \; \forall
n \in \omega \; (h(p_n) \leq m )$, then $\exists p \in \PP
\; \forall n \in \omega \; (p \leq p_n)$;
\par
\item{(II)} given $m \in \omega$ and $p,q \in \PP$ there is $\{
q_i ; \; i \in \ell \} \subseteq \PP$ so that
\par
\itemitem{(i)} $\forall i \in \ell \; (q_i \leq q \;\land\; q_i
\perp p )$;
\par
\itemitem{(ii)} whenever $q' \leq q$ is incompatible with $p$ and
$h(q') \leq m$ then there exists $i \in \ell$ so that
$q' \leq q_i$;
\par
\item{(III)} if $p,q \in \PP$ are compatible, there is $r \leq p,q$
so that $h(r) \leq h(p) + h(q)$.
\par
\noindent (Note that (by [BrJ, 1.7]) this notion is a strengthening
of the one given in [BrJ, 1.1].) We say a p.o. $\PP$ is
{\it soft} iff there is $\PP ' \subseteq r.o.(\PP)$ dense and
$h : \PP ' \to \omega$ so that $(\PP ' , h )$ is soft.
Furthermore, a pair $(\PP , h)$ {\it satisfies property $(*)$}
iff $\PP$ is a p.o., and $h$ is a height function on $\PP$
satisfying (III) above and:
\par
\item{$(*)$} given $p \in \PP$, a maximal antichain
$\{ p_n ; \; n \in \omega \} \subseteq \PP$ of conditions below
$p$ and $m \in \omega$, there exists $n \in \omega$ such
that: whenever $q \leq p$ is incompatible with $\{ p_j ;\;
j \in n\}$ then $h(q) > m$.
\par
\noindent A compactness argument shows:
\smallskip
{\capit Lemma 1.} (cf [BrJ, 1.2]) {\it If $(\PP , h)$ is soft,
then $(\PP , h)$ has property $(*)$.}
\smallskip
{\it Proof.} Put together the arguments of 1.7 and 1.2 in
[BrJ]. (This uses only (I) and (II) of softness.)
$\qed$
\smallskip
\noindent Our strategy to finish the proof of the Theorem
is as follows: show that $\PP_x$ is soft
(Lemma 4) --- and prove that iterating p.o.'s with
property $(*)$ doesn't increase ${\bf b}$ and ${\bf e}$
(Lemmata 2 and 3).
For the latter we need the following notion.
\par
A pair $\pi = ((A_k ;\; k\in\omega ) , (\pi_k ;\; k\in\omega ))$
is called a {\it generalized predictor} iff for all $k \in\omega$,
$A_k \subseteq [k , \omega )$ is finite, and $\pi_k$ is a 
function with $dom (\pi_k) = \{ \sigma \in \omega^{<\omega} ; \;
lh (\sigma) \in A_k \}$ and $ran (\pi_k) \subseteq
[\omega]^{<\omega}$. $\pi$ {\it predicts} $f \in \omega^\omega$
iff $\forall ^\infty k \; \exists \ell \in A_k \;
(f(\ell) \in \pi_k (f \restrict \ell))$; otherwise
$f$ {\it evades} $\pi$. --- The original definition of
predicting is a special instance of this notion
(in case $A_k = \{ \ell_k \}$, $\ell_k < \ell_{k+1}$
and $\vert \pi_k (\sigma) \vert =1$ for $\sigma
\in \omega^{\ell_k}$).
\par
Next let us consider the following property of p.o.'s $\PP$:
\smallskip
\item{(++)} given $F \subseteq \omega^\omega \cap V$, $F 
\in V$, such that for any countable family of generalized predictors
$\Pi$ there is $f\in F$ evading all $\pi \in\Pi$, and
$\langle \breve\pi_n ; \; n \in \omega \rangle$ a sequence
of $\PP$--names for generalized predictors, we can find a sequence
$\langle \pi_n ; \; n\in\omega\rangle \in V$ of generalized
predictors such that whenever $f \in F$ evades all
$\pi_n$, then
\par
\centerline{$\forces_\PP "f$ evades all $\breve\pi_n$".}
\smallskip
\noindent (Note that this is almost the same as $(**)$ in 2.1.)
This more general version of predicting as well as (++)
are needed for the following preservation results:
\smallskip
{\capit Lemma 2.} {\it Suppose $\PP$ is a $ccc$ p.o., $h$ is
a height function on $\PP$, and $(\PP , h)$ satisfies property
$(*)$. Then:
\par
\item{(a)} any unbounded family of functions in $\omega^\omega
\cap V$ is still unbounded in $V[G]$, where $G$ is
$\PP$--generic over $V$;
\par
\item{(b)} 
any family of functions in $\omega^\omega \cap V$ which
is not predicted by a single (countable family of)
generalized predictor(s) still has this property in $V[G]$ ---
and, in fact, $\PP$ satisfies $(++)$. \par}
\smallskip
{\capit Lemma 3.} {\it Assume $\langle \PP_\beta , \breve\QQ_\beta
; \; \beta < \alpha \rangle$ is an $\alpha$--stage finite
support iteration of $ccc$ p.o.'s such that $\forces_{\PP_\beta}
"h_\beta$ is a height function on $\breve\QQ_\beta"$.
Then:
\par
\item{(a)} if $\forall \beta < \alpha \; \forces_{\PP_\beta}
"\breve\QQ_\beta$ has property $(*)"$, then $\PP_\alpha$
does not add dominating reals. \par
\item{(b)} if $\forall \beta < \alpha \; \forces_{\PP_\beta}
"\breve\QQ_\beta$ has property $(*)"$ 
or just $\forces_{\PP_\beta} "\breve\QQ_\beta$ satisfies $(++)"$,
then $\PP_\alpha$ satisfies $(++)$ (and thus no
$\omega$-predictor
predicts all old reals in the Baire space). \par}
\smallskip
{\it Proof of Lemma 2.} (a) [BrJ, 1.3].
\par
{\openup2\jot
(b) Let $F$ be such a family in $\omega^\omega \cap V$. Suppose
$\forces_\PP " \breve\pi = ((\breve A_k ;\; k\in\omega ) ,(
\breve\pi_k ;\; k \in\omega ))$ is a generalized predictor".
For $k \in \omega$ let $\{ p^n_k ;\; n \in \omega \}$ be a
maximal antichain deciding the set $\breve A_k$. Choose $n_k$
according to $(*)$ so that: whenever $p$ is incompatible with 
$\{ p^j_k ; \; j \in n_k \}$, then $h(p) > k$.
For $j \in n_k$ let $A^j_k$ be such that $p^j_k \forces_\PP
"\breve A_k = A^j_k "$ and let $A_k := \bigcup_{j \in n_k}
A^j_k$.
\par
Fix $k \in \omega$ and $j \in n_k$. Let $A^j_k = \{
\ell^j_0 , ... , \ell^j_{a(j,k)-1} \}$ be the
increasing enumeration of $A^j_k$. By recursion on
$m < a(j,k)$ we define conditions $p^{j,\tau}_{k,\sigma}$
for $\sigma \in \omega^{\ell^j_m}$ and $\tau \in \omega^{m+1}$ as
follows:
in case $m = 0$, $\sigma \in \omega^{\ell^j_0}$, let
$\{ p^{j , \langle n \rangle }_{k , \sigma} ; \; n \in \omega
\}$ be a maximal antichain of conditions below $p^j_k$ deciding
the set $\breve\pi_k (\sigma)$; in case $m \geq 0$,
$\sigma \in \omega^{\ell^j_{m+1}}$, $\tau \in \omega^{m+1}$, let
$\{ p^{j,\tau\hat{\;}\langle n \rangle}_{k,\sigma} ; \;
n \in \omega \}$ be a maximal antichain of conditions below
$p_{k,\sigma\restrict\ell^j_m}^{j,\tau}$ deciding the
set $\breve\pi_k (\sigma)$. Next, for $\sigma \in \omega^{\ell^j_m}$,
$\tau \in \omega^m$, $m < a(j,k)$, choose $n^{j,\tau}_{k,\sigma}
\in \omega$ according to $(*)$ so that:
whenever $p < p^{j,\tau}_{k,\sigma\restrict \ell^j_{m-1}}$ 
($p<p^j_k$ in case $m=0$)
is incompatible with $\{ p^{j , \tau\hat{\;}\langle
i\rangle}_{k,\sigma} ; \; i \in n^{j,\tau}_{k,\sigma} \}$,
then $h(p) > k + h(p^j_k) + ... + h(p^{j,\tau}_{k,\sigma \restrict
\ell^j_{m-1}} )$.
For $i \in n^{j,\tau}_{k,\sigma}$ let $\pi^{j,\tau\hat{\;}
\langle i\rangle}_{k,\sigma}$ be such that $p^{j,\tau\hat{\;}
\langle i \rangle}_{k,\sigma} \forces_\PP" \breve\pi_k
(\sigma) = \pi^{j,\tau\hat{\;}\langle i \rangle}_{k,\sigma}"$.
Unfixing $j$, we define, for $\sigma \in \omega^{<\omega}$
with $dom(\sigma) \in A_k$,
$\pi_k (\sigma) : = \bigcup \{ \pi^{j,\tau}_{k,\sigma} ;
\; j \in n_k$ is such that $dom (\sigma) \in A^j_k$, e.g.,
$dom(\sigma) = \ell^j_m$ where $m < a(j,k)$, $\tau\in
\omega^{m+1}$, and for all
$i \leq m$, $\tau (i) \in n^{j, \tau \restrict i}_{k , \sigma
\restrict \ell^j_i} \}$. Then $\pi = ((A_k ;\; k \in\omega ),(
\pi_k ;\; k\in\omega ))$ is a generalized predictor.
Choose $f\in F$ evading $\pi$.
\smallskip
{\it Claim. $\forces_\PP "f$ evades $\breve\pi "$.}
\smallskip
{\it Proof of Claim.} Suppose there is a $p \in \PP$ and a $k_0
\in \omega$ so that
\smallskip
\centerline{$p \forces_\PP " \forall k \geq k_0 \; \exists \ell
\in \breve A_k \; (f(\ell) \in \breve\pi_k (f \restrict 
\ell ))"$.}
\smallskip
\noindent Choose $k \geq k_0$ so that $h(p) \leq k$ and
$\forall \ell \in A_{k} \; (f(\ell) \not\in \pi_k(f\restrict \ell
))$. Then $p$ is compatible with $p^j_{k}$ for some $j \leq
n_{k}$; let $q^j_{k}$ be a common extension such that $h
(q^j_{k}) \leq h(p) + h(p^j_{k} ) \leq k +
h(p^j_{k})$. Next construct recursively $\tau_m \in \omega^{m+1}
$ and $q^{j, \tau_m}_{k, f \restrict \ell^j_m}$ for $m < a
(j,k)$ as follows: in case $m = 0$, find
$i < n^j_{k,f\restrict\ell^j_0}$ so that $q^j_{k}$ and 
$p^{j,\langle i \rangle}_{k,f \restrict \ell^j_0}$ are compatible,
let $q^{j,\langle i \rangle}_{k,f\restrict\ell^j_0}$ be
a common extension of height $\leq k + h(p^j_{k})
+h(p^{j,\langle i\rangle}_{k, f\restrict \ell^j_0} ) $
and let $\tau_0 = \langle i\rangle$; in case $m \geq 0$,
find $i < n^{j,\tau_m}_{k,f\restrict\ell^j_{m+1}}$ so that
$q^{j,\tau_m}_{k,f \restrict \ell^j_m}$ and $p^{j,\tau_m
\hat{\;}\langle i\rangle}_{k,f\restrict \ell^j_{m+1}}$ are
compatible, let $q^{j,\tau_m\hat{\;}\langle
i\rangle}_{k,f\restrict
\ell^j_{m+1}}$ be a common extension of height $\leq
k + ... + h(p^{j,\tau_m\hat{\;}\langle i \rangle}_{k,
f \restrict \ell^j_{m+1}})$ and let $\tau_{m+1}
=\tau_m \hat{\;} \langle i \rangle$. Note that for
$m= a (j,k) - 1$ we have
\smallskip
\centerline{$q^{j,\tau_m}_{k,f \restrict \ell^j_m} \forces_\PP
"\forall \ell \in \breve A_k \; (f(\ell) \not\in \breve\pi_k
(f\restrict \ell))"$,}
\smallskip
\noindent a contradiction. $\qed$}
\smallskip
Note that a trivial modification of this argument proves
the
stronger version of Lemma 2 (b) as well. Notice also that
we used property (III) of the definition of softness in the proof
of the claim. $\qed$
\smallskip
{\it Proof of Lemma 3.} (a) [JS 2, 2.2]; see also [BrJ, 1.8].
\par
(b) Simply rewrite the proof of Lemma 2.2 in the present context.
$\qed$
\smallskip
Using these two preservation results as well as Lemma 1, we can
finish the proof of the Theorem by showing:
\smallskip
{\capit Lemma 4.} {\it For any $x \in \omega^\omega$,
$\PP_x$ is soft.}
\smallskip
{\it Proof.} Let $h: \PP_x \to \omega$ be defined by
$h(d, \langle \pi_k ; \; k\in d \rangle , F) :=
max \{ max (d) , \vert F \vert \}$. $h$ is trivially a
height function; furthermore any decreasing sequence of
$\PP_x$ which becomes eventually constant in height becomes
eventually constant in the first two coordinates, and in the third
coordinate, the size of the set of functions considered is
eventually constant --- so (I) of the definition of softness
is obvious. (III) is easy, and we are left with (II).
\par
We can assume $p = (d,\langle\pi_k ;\; k\in d \rangle
, F)$ and $q = (e , \langle \tilde\pi_k ;\; k\in e
\rangle ,G)$ are compatible and $q \not\leq p$.
We now describe which conditions of height $\leq m$ we put
into our finite set.
\par
(i) Assume $d \subseteq e$. Then we must have: $\exists
f \in F \; \forall g \in G \; (f \not\subseteq g)$.
We take all conditions of the form $(e' , \langle \pi_k ' ;\;
k \in e' \rangle , G )$ extending $q$ such that
$max (e') \leq m$ and for some $f \in F$ (with
$\forall g \in G \; (f \not\subseteq g)$) and
some $k \in e' \setminus e \; (\pi_k ' (f \restrict k ) \neq
f (k) )$.
\par
(ii) Assume $e \subset d$. We take all conditions of
the form $(e' , \langle \pi_k ' ;\; k \in e' \rangle , G')$
extending $q$ such that $max (e') \leq m$, $\vert G ' \vert \leq
m$ and either \par
\item{(a)} $G = G'$ and [$d \cap  max (e') 
\not\subseteq e'$ or $e' \cap  max (d) 
\not\subseteq d$] or \par
\item{(b)} $G = G'$ and for some $k \in (d \cap e')
\setminus e \; (\pi_k ' \neq \pi_k )$ or \par
\item{(c)} $G' = G \cup \tilde G$ where 
$\forall g \in\tilde G \; (dom(g) \leq max (d) + 1 \;\land\;
g \in \prod_{n<dom(g)} x(n))$
and $\exists g \in \tilde G
\; \exists k \in d \; (\pi_k (g \restrict k ) \neq g(k))$ or
\par
\item{(d)} $G = G'$ and $d \subset e'$ and for some $f \in F
$ (with $\forall g \in G \; (f \not\subseteq g)$) 
and some $ k \in e ' \setminus d \; (\pi_k ' (f \restrict k) \neq
f(k))$ [this is like (i)].
\par
\noindent In either case we have described a finite set of
conditions and leave it to the reader to check that each
condition of height $\leq m$ below $q$
incompatible with $p$ is
indeed below one of the conditions exhibited.
This finishes the proof of the Lemma. $\qed$
\smallskip
As in 2.2. we notice that we can in fact prove the
consistency of ${\bf e}_{ubd} = \lambda$ and ${\bf b}={\bf e}
=\kappa$ for arbitrary uncountable regular $\kappa < \lambda$.
\bigskip
{\sanse 3.4.} 
There is an even stronger result:
\smallskip
{\capit Theorem 1.} {\it It is consistent that ${\bf e}_{ubd}
=2^\omega =\omega_2$ and ${\bf d}=\omega_1$.}
\smallskip
\noindent To appreciate this recall that by [Bl 2, Theorem 13],
we have ${\bf e}\leq{\bf d}$. Our reason for nevertheless keeping
the result in 3.3. is that it allows us to choose the size of
${\bf e}, {\bf b}$ and ${\bf e}_{ubd}$ arbitrarily. The proof of
Theorem 1 goes as follows. Let
\smallskip
\centerline{$\lambda^* := min \{ \vert X \vert ; \;
\exists f \in (\omega \setminus 2)^\omega \; (X \subseteq
\prod_n f(n) \;\land\; \forall \phi \in \prod_n [f(n)]^n
\;\exists x\in X \;\exists^\infty n \; (x(n) \not\in \phi(n))) \}
.$}
\smallskip
\noindent Pawlikowski [Pa, Lemma 2.2.] proved that $\lambda^*$
equals the transitive additivity of the ideal ${\cal L}$.
Rewriting Blass' proof of $add({\cal L}) \leq {\bf e}$
[Bl 2, Theorem 12] in our context, one easily sees
$\lambda^* \leq {\bf e}_{ubd}$. Hence Theorem 1 follows
from Shelah's 
\smallskip
{\capit Theorem 2.} [Sh 326, section 2] {\it $CON (
\lambda^* = 2^\omega = \omega_2 \;\land\; {\bf d}
= \omega_1)$.}
\smallskip
\noindent (Note that this uses only one of the forcing notions
of Shelah's result, namely the one from [Sh 326, Proposition
2.9]. Also notice that Pawlikowski's p.o. [Pa, Theorem
2.4] is soft and thus his model for $\lambda^*
> {\bf b}$ is an alternative to ours for
showing Theorem 3.3. --- this wouldn't shorten the argument,
however, for we still would have to prove that iterating
soft p.o.'s doesn't increase {\bf e}.)
\bigskip
{\sanse 3.5.} {\it Evasion ideals and duality.}
With the concepts of evading and predicting we can associate
$\sigma$--ideals on the reals as follows.
Fix a space $X = \prod_n f(n)$, where $f \in (\omega+1 \setminus
2)^\omega$. Let ${\cal I}_X =\{ A \subseteq X;$ there
is a countable set of $X$--predictors $\Pi$ such that for
all $g \in A$ there is $\pi\in\Pi$ predicting $g \}$.
Note that the ideals ${\cal I}_X$ are subideals of the 
ideal ${\cal E}$ (see 3.1.). We shall study the
cardinal coefficients related to these ideals (see
$\S$ 1).
\smallskip
{\capit Proposition.} {\it $add({\cal I}_X) =
\omega_1$, $cof({\cal I}_X) = 2^\omega$ and $unif
({\cal I}_X)\geq {\bf e}_X$.}
\smallskip
{\it Proof.} Let $\{ D_\alpha ; \; \alpha <
2^\omega \}$ be an a.d. family of subsets of $\omega$.
Let $\pi_\alpha$ be an arbitrary $X$--predictor with
$D_{\pi_\alpha} = D_\alpha$. Let $A_\alpha =
\{ g \in X ;\; \pi_\alpha$ predicts $g\}$, and
note that whenever $\Gamma \subseteq 2^\omega$ is
uncountable then $\bigcup_{\alpha \in\Gamma}
A_\alpha \not\in {\cal I}_X$. This shows the
two equalities; the inequality is trivial.
$\qed$
\smallskip
\noindent It is unclear whether $unif({\cal I}_X) =
{\bf e}_X$ (cf section 5, question (4)); we note,
however, that all known $ZFC$--results as well as all
known consistency results about ${\bf e}_X$
carry over to $unif({\cal I}_X)$.
Also, the cardinal ${\bf e'}$ defined in 2.3. can
be viewed as the uniformity of a (slightly larger)
ideal. 
\par Similarly, dualizing these results (cf 1.1.), we
get corresponding results about $cov({\cal I}_X)$, the
smallest size of a set of predictors needed to predict all
reals from $X$. The only problem occurs when dualizing
the consistency results gotten from countable support iterations.
Concerning this we note that the Mathias real model satisfies
${\bf b} = \omega_2$ and $cov({\cal I}_X) = \omega_1$
for any space of the form $X=\prod_n f(n)$,
where $f \in (\omega\setminus 2)^\omega$ (this forms part of
the proof of Theorem D, see 4.3.), and thus is dual to
the one of 3.4., Theorem 1. On the other hand, to show the
consistency of $cov({\cal I}_{n^\omega}) = \omega_1$
and $cov({\cal I}_X) = unif ({\cal M}) = \omega_2$ for
$X = \prod_n f(n)$, $f\in (\omega\setminus 2)^\omega$
converging (fast enough) to infinity (dual to the
Mathias real model, see the remark at the end of subsection
3.2.), use the model gotten by
iterating with countable support the
forcing of [BaJS, section 3] (we leave the details of this to 
the reader; just note that the generic real is not predicted by
$X$--predictors from the ground model and thus $cov({\cal I}_X)
=\omega_2$ after the iteration, and that $cov({\cal I}_{n^\omega})
= \omega_1$ can be shown by an argument similar to the one that
the iteration doesn't add random reals [BaJS, 3.6, 3.8 and
3.15]).
\par
The relationship between the cardinals associated with
evading and predicting as well as many other cardinals can be
displayed in the following diagram (where the invariants
grow larger when moving up along the lines).
\smallskip
\centerline{put the diagram here}
\smallskip
\noindent Here ${\cal I} = {\cal I}_{\omega^\omega}$,
and ${\cal I}_\ell$ is the ideal associated with linear
predicting as defined by Blass; both ${\cal I}_\ell$
and ${\bf e}_\ell$ could be replaced by ${\cal I}_\KK$
and ${\bf e}_\KK$, respectively, where $\KK$ is a
countable field. To ease the reading we did not include
several inequalities not related to evading and predicting;
these are $add({\cal M}) \leq unif({\cal E})$ and 
$cov({\cal E}) \leq cof({\cal M})$, ${\bf b}\leq
unif({\cal M})$ and $cov({\cal M}) \leq {\bf d}$,
$cov({\cal L}) \leq unif({\cal M})$ and $cov({\cal M})
\leq unif ({\cal L})$, as well as ${\bf b}\leq{\bf r}$
and ${\bf s} \leq {\bf d}$. All inequalities are proved
in [vD, section 3], [Fr], [BS], [Bl 2] or our
work (see also [Va] for other references). 
Almost all inequalities are consistently strict
(see [vD, section 5], [BaJS] or our work); this is unclear
only for ${\bf e}_\ell \leq {\bf e}$ and dually
$cov({\cal I}) \leq cov({\cal I}_\ell)$; see question (2)
in section 5.
\Bigskip
{\dunhg $\S$ 4. Luzin sets of evading functions, Luzin groups,
and Gross spaces}
\Smallskip
{\sanse 4.1.} Recall (cf, e.g., [Mi 2, p. 206]) that given a
$\sigma$--ideal ${\cal I}$ on the real line, 
an {\it ${\cal I}$--Luzin set} is defined to be an uncountable subset of
the reals with at most countable intersection with
every member of ${\cal I}$. One of the goals of this
section is to study this notion in case of the
ideals introduced in 3.5.
\par
Let $X = \prod_n f(n)$, where $f \in (\omega + 1 \setminus
2 )^\omega$,
be one of the spaces studied in the latter subsection.
We say an uncountable $F \subseteq X$ is a {\it Luzin
set of evading functions} iff $F$ is ${\cal I}_X$--Luzin
iff for all $X$--predictors $\pi$ at most countably many $f \in F$
are predicted by $\pi$. More generally, given $\kappa$ with
$cf(\kappa) > \omega$, $F \subseteq X$ is a {\it generalized
Luzin set of evading functions of size $\kappa$}
iff $F$ is generalized ${\cal I}_X$--Luzin iff
for all $X$--predictors $\pi$ less that $\kappa$ many $f \in
F$ are predicted by $\pi$.
\par
It turns out, however, that in case of linear prediction
the following notion is more useful. Let $\KK$ be a finite
or countable field. An additive group $G = (G,+) \leq
(\KK^\omega , +)$ which is closed under multiplication
with elements from $\KK$ (i.e., it is a subspace of the
vector space $\KK^\omega$) is called a {\it Luzin group}
iff for all $\KK$--linear predictors $\pi$ at most
countably many $g \in G$ are predicted by $\pi$. Similarly, for 
$\kappa$ with $cf(\kappa) > \omega$, we say $G \leq \KK^\omega$
is a {\it generalized Luzin group of size $\kappa$} 
iff for all $\KK$--linear
predictors $\pi$ less than $\kappa$ many $g \in G$ are predicted
by $\pi$.
\bigskip
{\sanse 4.2.} Obviously, the existence of a Luzin set of evading
functions for a space $X$ implies that the corresponding
evasion number equals $\omega_1$. Using ideas from [JS 3],
we proceed to show that the converse need not hold.
\smallskip
{\capit Theorem.} {\it It is consistent that ${\bf e}_{fin} =
\omega_1$, but there are no Luzin sets of evading functions
for any space $n^\omega$ (where $n \in \omega$) and no
Luzin groups $G \leq \KK^\omega$, where $\KK$ is the two--element
field.
Actually, the conclusion holds in Laver's model
for the Borel conjecture.}
\smallskip
{\it Proof.} ${\bf e}_{fin} = \omega_1$ follows from ${\bf e}_{fin}
\leq unif ({\cal L})$ and the fact that $unif({\cal L}) =
\omega_1$ holds in Laver's model [JS 2, section 1].
\par
Using a similar argument as in 3.2. (Lemma 2) it is easy to see that
the existence of a Luzin set of evading functions for $n^\omega$
($n \geq 2$) is equivalent to the existence of a Luzin set of
evading functions for the Cantor space $2^\omega$. Furthermore,
as any such Luzin set of size $
\omega_1$ would lie in an intermediate extension
of Laver's model, it suffices to show that adding one Laver real
destroys the ground model's Luzin sets.
\par
So assume $F \in V$ is a Luzin set of evading functions
(for $2^\omega$). We introduce an $\LL$--name $\breve \pi$
for a predictor as follows: if $\breve \ell$ is the $\LL$--name
for the generic real, let $\breve D_{\breve \pi}
:= \{ \breve \ell (i) ; \; i \in \omega \}$ be the name for the
set of numbers on which we predict, and let $\breve
\pi_k (\sigma) = 0 $ for all $\sigma \in 2^k$ and $k \in
\breve D_{\breve \pi}$. 
\smallskip
{\it Claim. $\forces_\LL$ "$\breve \pi$ predicts uncountably
many elements from $F$".}
\smallskip
{\it Proof of Claim.} Let $T \in \LL$, $N \prec \langle H(\kappa)
, ... \rangle $ countable with $T \in N$. Let $f \in F$ be
such that $f$ evades all predictors from $N$ (note that all but
countably many elements of $F$ have this property). We shall
construct $S \leq T$, $S \in \LL$ such that
\smallskip
\centerline{$S \forces_\LL $"$\breve \pi$ predicts $f$".}
\smallskip
\noindent Clearly this is enough (it shows that no $T \in \LL$
forces that $\breve \pi$ predicts only countably many
elements).
\par
For $\rho \in T$ such that $succ_T(\rho)$ is infinite, let
\smallskip
\centerline{$A_\rho := \{ n \in succ_T (\rho ) ; \; f(n) =
0 \}$.}
\smallskip
\noindent Note that $A_\rho$ must be infinite (otherwise,
define a predictor $\pi$ by: $D_\pi : = succ_T(\rho)$ and
$\pi_k (\sigma) = 1$ for $\sigma \in 2^k$ 
($k \in D_\pi$) iff $\sigma(max(D_\pi \cap k))=1$;
then $\pi \in N$ and $\pi $ predicts
$f$, a contradiction).
Now define $S$ by recursion on its levels:
\par
\itemitem{(i)} $stem(S) = stem(T)$; \par
\itemitem{(ii)} assume $\rho \in S$, then:
$\rho \hat{\;} \langle n \rangle \in S \longleftrightarrow
n \in A_\rho$.
\par
\noindent Note that $S \forces_\LL " \breve \pi$ predicts
$f$". This concludes the proof of the claim. $\qed$
\smallskip
To see that there are no Luzin groups over $\KK$, note that both
the predictor defined from the Laver real in the extension
and the predictor
defined in the proof of the claim are in fact linear. 
$\qed$
\bigskip
{\sanse 4.3.} We start the proof of the dichotomy theorem
(Theorem D) and recall (an extended version of)
the statement of its first part.
\smallskip
{\capit Theorem.} {\it It is consistent that
$2^\omega = \omega_2$ and there are no generalized Luzin sets
of evading functions for any space $n^\omega$ (where
$2 \leq n \in \omega$). Furthermore there are no generalized
Luzin groups $G \leq \KK^\omega$, where $\KK$ is a finite field.}
\smallskip
{\it Proof.} This is true in the Mathias
model.
That there are no Luzin sets (Luzin groups) of size
$\omega_1$ follows 
from ${\bf e}_{fin} \geq {\bf e_\KK} \geq {\bf s}
$ (3.2., Lemma 3) and ${\bf s} = \omega_2$.
\par
To see that there are no generalized Luzin sets (Luzin groups) of
size $\omega_2$ we apply the Laver property (see 
1.2.). In fact, we show
something slightly more general:
the ground model predictors predict all new reals
of the space $X=\prod_n f(n)$, where $f\in\omega^\omega$ is
arbitrary
(i.e. $cov({\cal I}_X) = \omega_1$ in the language of
3.5.).
Thus, any set of size $\omega_2$ contains a subset of same
size predicted by one predictor. 
\par
Let $\PP$ be a proper forcing notion satisfying the Laver
property. Let $\breve g$ be a $\PP$--name for an element
of $X$. Partition $\omega$ into disjoint intervals $\{
I_n ; \; n \in \omega \}$ such that $\vert I_n \vert = n^2$,
and $max (I_n) + 1 = min (I_{n+1})$. Using the Laver property,
we can find for given $p \in \PP$ a $q \leq p$ and a sequence 
$\{ \sigma^k_n ; \; n \in \omega \land k \in n \}$ such 
that $$\forall n \in \omega \; \forall k \in n \;
(\sigma^k_n \in \prod_{m \in I_n} f(m))$$
and 
\smallskip
\centerline{$q \forces_\PP " \forall n \; \exists k \in
n \; (\breve g \restrict I_n = \sigma^k_n)"$.}
\smallskip
\noindent Fix $n$. As $\vert I_n \vert = n^2$ there must be
$m_n \in I_n$ such that for $k_1 , k_2 \in n$:
whenever $\sigma^{k_1}_n \restrict m_n = \sigma_n^{k_2} \restrict
m_n$, then $\sigma^{k_1}_n (m_n) = \sigma^{k_2}_n (m_n)$.
Let $D_\pi = \{m_n ; \; n \in \omega\}$ and define the predictor
$\pi$ by: for $n \in \omega$ and
$\sigma \in \prod_{i \in m_n} f(i)$,
$$\pi_{m_n}(\sigma) = \cases{ \sigma_n^k (m_n) & if
$\sigma \restrict (I_n \cap m_n ) =\sigma^k_n \restrict m_n$
\cr 0 & otherwise.\cr} \leqno(*)$$
By the choice of $m_n$, $\pi_n$ is well--defined. It follows from
the construction that 
\smallskip
\centerline{$q \forces_\PP " \pi$ predicts $\breve g$".}
\smallskip
Finally, in case of generalized Luzin groups, we have to define
a linear predictor in the ground model. To do this, we proceed
as above, and take for fixed $n$ a maximal linearly independent
subset of the $\sigma^k_n$, $k \in n$, and use this in the
definition
corresponding to $(*)$. $\qed$
\bigskip
{\sanse 4.4.} {\it End of proof of dichotomy theorem.}
We got the idea to prove part (b) of Theorem D from [ShSp,
section 3].
\par
Let $\{ \phi_n ; \; n \in \omega \}$ enumerate
all linear functions from some $\KK^m$ to $\KK$
($m \in \omega$). Let $\KK = \{ a_n ; \; n \in \omega \}$.
Furthermore choose an unbounded, well--ordered
(under $\leq^*$) family $\{ f_\alpha ; \; \alpha < {\bf b}
\}$ of strictly increasing functions from $\omega$
to $\omega$.
\par
For $\sigma \in \omega^{<\omega}$ we define $\sigma^*
\in \omega^{<\omega}$ by recursion on $lh(\sigma)$ and
$\sigma(lh(\sigma) - 1)$ as follows. Fix $n , m \in \omega$
and $\sigma \in \omega^n$ such that $\sigma( n-1 ) = m$. For
$i < n-1 $ we let
\smallskip
\centerline{$\sigma^* (i) := (\sigma \restrict (i+1 ))^* (i)$.}
\smallskip
\noindent We let $\sigma^* (n-1)$ be the minimal $\ell$
such that
\smallskip
\centerline{(+) $a_\ell \not\in \{ \phi_k (a_{\ell_0} ,...);
k< m \land \ell_j \in \cup\{ran(\tau^*) ;
\tau \in m^{\leq n}\} ({\rm where} \; j < dim(dom(\phi_k)))\}$.}
\smallskip
\noindent This concludes the definition of $\sigma^*$.
\par
Next, given $\alpha < {\bf b}$, we define $g_\alpha \in
\KK^\omega$ by:
\smallskip
(++)\centerline{$g_\alpha (n) := a_{(f_\alpha \restrict (n+1))^*
(n)}.$}
\smallskip
\noindent We claim that $G := \langle g_\alpha ; \; \alpha
< {\bf b} \rangle$ is a Luzin group.
\par
For suppose not. Then (without loss) there are $n \in \omega$,
$A_\alpha \subseteq {\bf b}$ ($\alpha < {\bf b}$) a 
$\Delta$--system of sets of size $n$, $b_\ell \in
\KK$ ($\ell < n$) and a linear predictor $(D_\pi,
(\pi_m ; \; m \in D_\pi ))$ such that
$$\forall m \in D_\pi \; \forall \alpha < {\bf b}
\; ( \pi_m (\sum_{\ell < n} b_\ell g_{A_\alpha (\ell)}
\restrict m ) = \sum_{\ell < n} b_\ell g_{A_\alpha (\ell)}
(m)),$$
$$(*) \;\;\; {\rm i.e.}\; {1\over b_{n-1}} \cdot
\pi_m (\sum_{\ell < n} b_\ell g_{A_\alpha (\ell)} (0)
,..., \sum_{\ell < n} b_\ell g_{A_\alpha (\ell)}
(m-1)) - {1 \over b_{n-1}} \cdot \sum_{\ell
< n-1} b_\ell g_{A_\alpha (\ell)} (m) =
g_{A_\alpha (n-1)} (m).$$
Without loss we may as well assume that for some $m_0 \in
\omega$ for all $\alpha < {\bf b}$, for all $m \geq m_0$
and for all $\ell < n-1$, $f_{A_\alpha (n - 1)}
(m) > f_{A_\alpha (\ell)} (m)$.
Now choose $m \geq m_0$, $m \in D_\pi$ such that
$\{ f_{A_\alpha (n-1)} (m) ; \; \alpha < {\bf b} \}$ is
unbounded in $\omega$ (this is possible since $\{
f_{A_\alpha (n-1) } ; \; \alpha < {\bf b} \}$ is unbounded
in $\omega^\omega$, because $\{f_\alpha ; \; \alpha < {\bf b}
\}$ is unbounded in $\omega^\omega$ and well--ordered),  
and think of the expression $(*)$ as a linear function from
$\KK^{m \cdot n + n - 1}$ to $\KK$. Next choose $k \in \omega$
such that $\phi_k$ is the left--hand side
of $(*)$ and $\alpha < {\bf b}$ such that
$k < f_{A_\alpha (n-1)} (m)$. Then, for $\ell < n-1$,
$f_{A_\alpha (\ell)} \restrict (m+1) \in (f_{A_\alpha (n-1)}
(m) )^{m+1}$ and $f_{A_\alpha (n-1)} \restrict m
\in (f_{A_\alpha (n-1)} (m))^m$, and it is immediate 
from (+) and (++) that $g_{A_\alpha (n-1)} (m)
\neq \phi_k (g_{A_\alpha (0) }(0) , ... , g_{A_\alpha (n-2)}
(m) , g_{A_\alpha (n-1)} (0) , ... ,$ $g_{A_\alpha (n-1)}
(m-1))$, contradicting equation $(*)$. $\qed$
\smallskip
Note that the statement of the Theorem is a strong way of
saying ${\bf e}_\KK \leq {\bf b}$. For $\KK = \QQ$ this
was proved (rather indirectly using various 
intermediate group--theoretical notions) by Blass (see
Introduction, cf also 3.1.).\par
Granted the inequality ${\bf e}_\KK \leq {\bf b}$, there is,
in a sense, nothing peculiar about this
result. It merely reflects the fact that there is (in ZFC)
a set with Luzin--style properties associated with ${\bf
b}$ (this is the unbounded and well--ordered family
$\{f_\alpha; \; \alpha < {\bf b} \}$ we started with) and that
it is a (seemingly) general state of affairs that given such
a Luzin--style set for one cardinal (in our case, {\bf b})
we can construct a similar set (of same size) for a smaller
cardinal (in our case, ${\bf e}_\KK$) (see [Ci, Theorem 3.1.]
for related results). The main difficulty then is
to get a Luzin {\it group} and not just a Luzin set
of evading functions for the family of linear
$\KK$--valued predictors (cf $\S$ 5, question (5)).
\bigskip
{\sanse 4.5.} We prove a more general version of
the implication $(a) \Longrightarrow (b)$ in the
equivalence theorem (Theorem E).
\smallskip
{\capit Lemma.} {\it (A) If there is a 
strong Gross space $(E, \Phi)$ of dimension $\kappa$
over the field $\KK$, then there is a Luzin group
$G \leq \KK^\omega$ of size $\kappa$. \par
(B) Assume $cf(\kappa) > \omega$. 
If there is a Gross space $(E, \Phi)$
of dimension $\kappa$ over the field $\KK$, then there is
a generalized Luzin group $G \leq \KK^\omega$
of size $\kappa$.}
\smallskip
{\it Proof.} Both are similar. So we shall only prove (A) and
leave (B) to the reader. \par
Let $(E, \Phi)$ be strongly Gross of dimension $\kappa$
over $\KK$. Assume $\{ e_\alpha ; \; \alpha < \kappa \}$
is a basis of $\KK$. For $\alpha \geq \omega$ we define
$f_\alpha : \omega \to \KK$ by
\smallskip
\centerline{$f_\alpha (n) : = \Phi (e_n, e_\alpha )$.}
\smallskip
\noindent We claim that the subspace $G \leq \KK^\omega$
generated by the $f_\alpha$ ($\omega \leq \alpha < \kappa$)
is a Luzin group.
\par
For suppose not. Then there is a linear $\KK$--valued
predictor $\pi$ predicting $\omega_1$ many $g_\alpha
\in G$ ($\alpha < \omega_1$). Let $\pi = (D_\pi
, (\pi_m ; \; m \in D_\pi ))$; without loss, for all
$m \in D_\pi$ for all $\alpha < \omega_1$,
we have $\pi_m (g_\alpha \restrict m ) = g_\alpha (m)$.
We can assume that there are $n \in \omega$,
$A_\alpha \subseteq \kappa \setminus \omega$ ($\alpha < \omega_1$)
forming a $\Delta$--system of sets of size $n$,
and $b_\ell \in \KK$ ($\ell < n$) such that
\smallskip
\centerline{$g_\alpha = \sum_{\ell < n} b_\ell f_{A_\alpha
(\ell)}$ for $\alpha < \omega_1$.}
\smallskip
\noindent For $m \in D_\pi$ and $i \in m$, we let
$c_{im} := \pi_m (\sigma_{im}) \in \KK$ where
$\sigma_{im} \in \KK^m$ is such that
$$\sigma_{im} (j) = \cases{1 &if $i=j$ \cr
0 &otherwise. \cr}$$
By linearity of $\pi_m$, we have for any $\sigma \in \KK^{m+1}$
satisfying $\pi_m (\sigma \restrict m) = \sigma (m)$:
\smallskip
\centerline{$\sigma (m) = \sum_{j < m} c_{jm} \sigma (j)$,}
\smallskip
\noindent i.e.
\smallskip
\centerline{$\sum_{j < m} c_{jm} \sigma (j) - \sigma (m) = 0$.}
\smallskip
\noindent Thus we have for all $m \in D_\pi$ and all $\alpha
< \omega_1$
$$0 = \sum_{j < m} c_{jm} \cdot g_\alpha (j) - g_\alpha (m) =
\sum_{\ell < n} b_\ell (\sum_{j<m} c_{jm} \cdot f_{A_\alpha (\ell)}
(j) - f_{A_\alpha (\ell)} (m)) =$$
$$ = \sum_{\ell < n} b_\ell (\sum_{j<m} c_{jm} \cdot \Phi (e_j , 
e_{A_\alpha
(\ell)}) - \Phi (e_m , e_{A_\alpha (\ell)}) )=
\Phi (\sum_{j<m} c_{jm} \cdot e_j - e_m , \sum_{\ell < n}
b_\ell \cdot e_{A_\alpha (\ell)}).$$
Hence, if $U$ is the subspace of $E$ spanned by the vectors
$\sum_{j<m} c_{jm} \cdot e_j - e_m$, $m \in D_\pi$,
then $U^\perp \geq \langle \sum_{\ell < n} b_\ell \cdot
e_{A_\alpha (\ell)} ; \; \alpha < \omega_1 \rangle$.
Therefore $(E, \Phi)$ cannot be strongly Gross,
a contradiction. $\qed$
\bigskip
{\sanse 4.6.} {\it The proof of (b) $\Longrightarrow$
(a) in the equivalence theorem (Theorem E).}
The following argument was heavily influenced by
[ShSp, section 4, Theorem 4].
\par
Let $G \leq \KK^\omega$ be a Luzin group of size
$\omega_1$. Let $\{ g_\alpha ; \; \alpha < \omega_1 \}$
be a set of generators of $G$ as a $\KK$--vector space.
\par
Choose recursively injective functions $h_\alpha : \alpha
\to \omega$ such that $\omega \setminus ran(h_\alpha)$
is infinite and $\forall \beta < \alpha$ the set
$\{ \gamma < \beta ; \; h_\beta (\gamma) \neq 
h_\alpha (\gamma) \}$ is finite (this is one of the
standard constructions of an Aronszajn tree,
due to Todorcevic (cf [To, (2.2.)]))
\par
Let $E$ be a vector space of dimension $\omega_1$ over
$\KK$; assume that $E = \langle e_\alpha ; \; \alpha 
< \omega_1 \rangle$. We define a symmetric bilinear form
$\Phi$ on $E$ as follows:
\smallskip
\centerline{$\Phi (e_\alpha , e_\beta ) := g_\beta
(h_\beta (\alpha))$ for $\alpha < \beta <\omega_1$.}
\smallskip
\noindent We claim that $(E, \Phi)$ is a Gross space.
For if this were not the case, we could find (using
standard thinning--out arguments) vectors $y_k$
($k \in \omega$) and $z_\gamma$ ($\gamma \in \omega_1$),
and $\alpha^* \in \omega_1$, $n \in \omega$, 
$m_k$ ($k \in \omega$), $B_\gamma$ ($\gamma \in \omega_1$),
$A_k$ ($k \in \omega$), $b_i$ ($i \in n$) and $a_{jk}$
($j \in m_k$ and $k \in \omega$) such that
\smallskip
\itemitem{(i)} for all $k \in \omega$ and $\gamma \in
\omega_1$ we have $\Phi (y_k , z_\gamma) = 0$;
\par
\itemitem{(ii)} $A_k \subseteq \alpha^*$, $\vert A_k 
\vert = m_k$, and $a_{jk} \in \KK$ such that
$y_k = \sum_{j \in m_k} a_{jk} e_{A_k (j)}$; furthermore
$k_1 < k_2$ implies $max h_{\alpha^*} (A_{k_1})
< min h_{\alpha^*} (A_{k_2})$;
\par
\itemitem{(iii)} $B_\gamma \subseteq \omega_1$,
$\vert B_\gamma \vert = n$, and $b_i \in \KK$
such that $z_\gamma = \sum_{i\in n} b_i e_{B_\gamma
(i)}$; furthermore $\gamma_1 < \gamma_2 $ implies
$\alpha^* < min(B_{\gamma_1}) < 
max(B_{\gamma_1}) < min (B_{\gamma_2})$.
\smallskip
\noindent Next let us introduce a linear $\KK$--valued
predictor $\pi$ as follows. Fix $k \in \omega$.
Let $d_k := max h_{\alpha^*} (A_k)$; set
$D_\pi := \{d_k ; \; k \in \omega \}$; and let 
$j_k$ be such that $d_k = h_{\alpha^*} (A_k (j_k))$.
We set
$$\pi_k (\sigma_{ik}) := \cases{- {a_{jk} \over a_{j_k k}}
&if $i = h_{\alpha^*} (A_k (j))$ \cr
0 &if $i \not\in h_{\alpha^*} (A_k)$, \cr}$$
where $\sigma_{ik} \in \KK^{d_k}$ is such that
$$\sigma_{ik} (j) = \cases{1 &if $i=j$ \cr
0 &otherwise. \cr}$$
We extend $\pi_k$ linearly to $\KK^{d_k}$, and thus define a linear
$\KK$--valued predictor. By Luzinity of $G$, there is a
$\gamma \in \omega_1$ such that $g := \sum_{i\in n}
b_i g_{B_\gamma (i)}$ evades $\pi$. Choose $k \in \omega$
such that \smallskip
\itemitem{(1)} $h_{B_\gamma (i)} \restrict A_k = h_{\alpha^*}
\restrict A_k$ for all $i \in n$; \par
\itemitem{(2)} $\pi_k (g\restrict d_k) \neq g(d_k)$.
\smallskip
\noindent On the other hand, we have
$$0 = \Phi (y_k , z_\gamma) = \sum_{j \in m_k}
\sum_{i \in n} b_i \cdot a_{jk} \cdot g_{B_\gamma (i)}
(h_{B_\gamma (i)} (A_k (j))) =
\sum_{j \in m_k} \sum_{i \in n} b_i \cdot a_{jk}
\cdot g_{B_\gamma (i)} (h_{\alpha^*} (A_k (j))).$$
I.e.
$$\sum_{j \in m_k \setminus \{j_k\} } a_{jk}
\cdot g(h_{\alpha^*} (A_k (j))) = - a_{j_kk} g(d_k).$$
Thus
$$\pi_k (g \restrict d_k) = - \sum_{j \in m_k
\setminus \{ j_k \} } {a_{jk} \over a_{j_kk}}
g(h_{\alpha^*} (A_k (j))) = g(d_k),$$
a contradiction. $\qed$
\bigskip
{\sanse 4.7.} {\capit Corollary.} (Baumgartner, Shelah,
Spinas; [BSp], [ShSp], see also [Sp 2])
{\it Let $\KK$ be an arbitrary finite or countable field.
\par
\item{(a)} Assume any of the following: \par
\itemitem{} --- there is a Luzin group $G \leq \KK^\omega$
\par
\itemitem{} --- $cof ({\cal M}) = \omega_1$
\par
\itemitem{} --- ${\bf b} = \omega_1$ (in case $\vert \KK
\vert = \omega$)
\par
\item{} Then there is a strong Gross space over $\KK$.
\par
\item{(b)} Assume any of the following: \par
\itemitem{} --- there is no Luzin group $G \leq \KK^\omega$
\par
\itemitem{} --- ${\bf e}_\KK > \omega_1$ \par
\itemitem{} --- $add ({\cal L}) > \omega_1$ \par
\itemitem{} --- ${\bf p} > \omega_1$ \par
\itemitem{} --- ${\bf s} > \omega_1$ (in case $\vert
\KK \vert < \omega$) \par
\item{} Then there is no strong Gross space over $\KK$.
}
\smallskip
{\it Proof.} (a) 4.4. and 4.6. \par
\noindent
We leave the construction of a Luzin group from $cof ({\cal M})
= \omega_1$ to the reader.
\par
(b) 4.5. and [Bl, section 4] --- see also $\S$ 3 and in
particular 3.2. $\qed$
\smallskip
In a sense our results say that there is no cardinal
invariant such that its being $\omega_1$ is equivalent
to the existence of strong Gross spaces over $\KK$ (cf
the question in the Introduction). The natural candidate
for such a cardinal would be ${\bf e}_\KK$, but, by
4.2., ${\bf e}_\KK$ may be $\omega_1$ and there may be
no Luzin subgroup of $\KK^\omega$ and hence no
strong Gross space over $\KK$.
\par
We note in closing that we think of the dichotomy theorem
as the basic result underlying the fact that it is much more
difficult to get rid of Gross spaces over countable fields
than over finite fields (in fact, the consistency
of the non--existence of Gross spaces over countable
fields is an open problem [Sp 2]).
It follows from 4.5. and 4.3. that there are no
Gross spaces over finite fields in Laver's or Mathias'
models (this was known to be true previously 
[ShSp, section 4, Theorem 2] in a model constructed
by Shelah [BlSh] in which there are both $P_{\omega_1}$--
and $P_{\omega_2}$--points --- in this peculiar situation it
is indeed not difficult to see that there are no such
spaces; however we think that Laver's or Mathias' models
are much easier to grasp combinatorially).

\Bigskip
{\dunhg $\S$ 5. Questions}
\Smallskip
We have introduced a multitude of cardinals and in spite
of our $ZFC$-- and consistency results many questions
concerning the relationship between them remain open.
We mention but a few.
\par
The results of section 2 suggest a positive answer to the
following.
\smallskip
\item{(1)} {\it Is ${\bf se} \leq add({\cal M})$
(or even $min \{ {\bf e'} ,{\bf b} \} \leq add ({\cal M})$)
in $ZFC$? Is ${\bf se} < add({\cal M})$ (or even ${\bf e'}
< add({\cal M})$) consistent?}
\smallskip
\noindent Recall that Blass proved $min\{ {\bf e} , {\bf b}
\} \leq add ({\cal M})$ [Bl 2, Theorem 13]. To calculate the
value of ${\bf e'}$ in the Hechler real model may shed some light 
on the situation.
\par
The most important problem is perhaps:
\smallskip
\item{(2)} (Blass [Bl 2, section 5, question (2)]) {\it Clarify
the relationship between {\bf e} and {\bf b} (and between {\bf e}
and ${\bf e}_\ell$)! Or: is there a generalized Luzin set
of evading functions for $\omega^\omega$ (cf 4.3. and 4.4.)?}
\smallskip
\noindent
Related is
\smallskip
\item{(3)} {\it Clarify the relationship between the different
${\bf e}_\KK$'s (and between ${\bf e}_\KK$ and ${\bf e}_{fin}$
for finite $\KK$)!
Or: does the existence of a strong Gross space over some finite
field imply the existence of a strong Gross space over
every finite or countable field?}
\smallskip
Concerning sections 3 and 4, the following additional 
questions may be of some interest:
\smallskip
\item{(4)} {\it Is ${\bf e}_X = unif ({\cal I}_X)$ (cf 3.5.)?}
\smallskip
\item{(5)} {\it Does the existence of a Luzin or Sierpi\'nski
set imply the existence of a Luzin group
(cf 4.4.)?}
\Bigskip
\centerline{\capitg References}
\Smallskip
\itemitem{[BaJS]} {\capit T. Bartoszy\'nski, H. Judah and S.
Shelah,} {\it The Cicho\'n diagram,} submitted to Journal of Symbolic
Logic.
\smallskip
\itemitem{[BS]} {\capit T. Bartoszy\'nski and S. Shelah,}
{\it Closed measure zero sets,} Annals of Pure
and Applied Logic, vol. 58 (1992), pp. 93-110.
\smallskip
\itemitem{[Bau]} {\capit J. Baumgartner,} {\it Iterated forcing,}
Surveys in set theory (edited by A.R.D. Mathias), Cambridge
University Press, Cambridge, 1983, pp. 1-59.
\smallskip
\itemitem{[BD]} {\capit J. Baumgartner and P. Dordal,} {\it
Adjoining dominating functions,} {Journal of Symbolic Logic,} vol.
50 (1985), pp. 94-101.
\smallskip
\itemitem{[BSp]} {\capit J. Baumgartner and O. Spinas,}
{\it Independence and consistency proofs in quadratic form theory,}
Journal of Symbolic Logic, vol. 56 (1991), pp. 1195-1211.
\smallskip
\itemitem{[Bl 1]} {\capit A. Blass,} {\it Simple cardinal
characteristics of the continuum,} in: Set theory of the reals,
Proceedings of the Bar--Ilan conference in honour of
Abraham Fraenkel.
\smallskip
\itemitem{[Bl 2]} {\capit A. Blass,} {\it Cardinal characteristics
and the product of countably many infinite cyclic groups,}
preprint.
\smallskip
\itemitem{[BlSh]} {\capit A. Blass and S. Shelah,}
{\it There may be simple $P_{\aleph_1}$-- and $P_{\aleph_2}$--points
and the Rudin--Keisler ordering may be downward directed,}
Annals of Pure and Applied Logic, vol. 33 (1987),
pp. 213-243.
\smallskip
\itemitem{[Br]} {\capit J. Brendle,} {\it Larger cardinals
in Cicho\'n's diagram,} Journal of Symbolic Logic,
vol. 56 (1991), pp. 795-810.
\smallskip
\itemitem{[BrJ]} {\capit J. Brendle and H. Judah,} {\it Perfect
sets of random reals,} to appear in Israel Journal of Mathematics.
\smallskip
\itemitem{[BrJS]} {\capit J. Brendle, H. Judah and S. Shelah,}
{\it Combinatorial properties of Hechler forcing,} Annals
of Pure and Applied Logic, vol. 58 (1992), pp. 185-199.
\smallskip
\itemitem{[Ci]} {\capit J. Cicho\'n,} {\it On two--cardinal
properties of ideals,} Transactions of the American Mathematical
Society, vol. 314 (1989), pp. 693-708.
\smallskip
\itemitem{[E]} {\capit K. Eda,} {\it A note on subgroups
of $\ZZ^\NN$,} Abelian Group Theory, Proceedings, Honolulu
1982/83 (R. G\"obel, L. Lady and A. Mader, eds.),
Lecture Notes in Mathematics 1006, Springer--Verlag,
1983, pp. 371-374.
\smallskip
\itemitem{[Fr]} {\capit D. Fremlin,} {\it Cicho\'n's diagram,}
S\'eminaire Initiation \`a l'Analyse (G. Choquet,
M. Rogalski, J. Saint Raymond), Publications Math\'ematiques
de l'Universit\'e Pierre et Marie Curie, Paris, 1984,
pp. 5-01 - 5-13. 
\smallskip
\itemitem{[Je 1]} {\capit T. Jech,} {\it Set theory,} Academic Press,
San Diego, 1978.
\smallskip
\itemitem{[Je 2]} {\capit T. Jech,} {\it Multiple forcing,}
Cambridge University Press, Cambridge, 1986.
\smallskip
\itemitem{[JS 1]} {\capit H. Judah and S. Shelah,} {\it
Souslin forcing,} Journal of Symbolic Logic, vol. 53 (1988),
pp. 1188-1207.
\smallskip
\itemitem{[JS 2]} {\capit H. Judah and S. Shelah,} {\it The Kunen-Miller
chart (Lebesgue measure, the Baire property, Laver reals and
preservation theorems for forcing),} Journal of
Symbolic Logic, vol. 55 (1990), pp. 909-927.
\smallskip
\itemitem{[JS 3]} {\capit H. Judah and S. Shelah,}
{\it Killing Luzin and Sierpi\'nski sets,} to appear.
\smallskip
\itemitem{[Ku]} {\capit K. Kunen,} {\it Set theory,} North-Holland,
Amsterdam, 1980.
\smallskip
\itemitem{[Mi 1]} {\capit A. Miller,} {\it Some properties of
measure and category,} Transactions of the American Mathematical
Society, vol. 266 (1981), pp. 93-114.
\smallskip
\itemitem{[Mi 2]} {\capit A. Miller,} {\it Special subsets of the
real line,} Handbook of set--theoretic topology, K. Kunen and 
J. E. Vaughan (editors), North--Holland, Amsterdam, 1984, pp.
201-233.
\smallskip
\itemitem{[Pa]} {\capit J. Pawlikowski,} {\it Powers of
transitive bases of measure and category,} Proceedings
of the American Mathematical Society, vol. 93 (1985),
pp. 719-729.
\smallskip
\itemitem{[Sh 326]} {\capit S. Shelah,} {\it Vive la diff\'erence
I: nonisomorphism of ultrapowers of countable models,}
Set Theory of the Continuum (H. Judah, W. Just and H. Woodin,
eds.), Mathematical Sciences Research Institute Publications,
Springer--Verlag, 1992, pp. 357-405.
\smallskip
\itemitem{[ShSp]} {\capit S. Shelah and O. Spinas,}
{\it How large orthogonal complements are there in a quadratic
space?} preprint.
\smallskip
\itemitem{[S]} {\capit E. Specker,} {\it Additive Gruppen von
Folgen ganzer Zahlen,} Portugaliae Mathematica, vol. 9 (1950),
pp. 131-140.
\smallskip
\itemitem{[Sp 1]} {\capit O. Spinas,} {\it Iterated forcing in
quadratic form theory,} Israel Journal of Mathematics,
vol. 79 (1992), pp. 297-315.
\smallskip
\itemitem{[Sp 2]} {\capit O. Spinas,} {\it Cardinal invariants
and quadratic forms,} in: Set theory of the reals, Proceedings
of the Bar--Ilan conference in honour of Abraham Fraenkel.
\smallskip
\itemitem{[To]} {\capit S. Todorcevic,} {\it Partitioning
pairs of countable ordinals,} Acta mathematica, vol. 159
(1987), pp. 261-294.
\smallskip
\itemitem{[vD]} {\capit E. K. van Douwen,} {\it
The integers and topology,} Handbook of set--theoretic
topology, K. Kunen and J. E. Vaughan (editors),
North--Holland, Amsterdam, 1984, pp. 111-167.
\smallskip
\itemitem{[Va]} {\capit J. Vaughan,} {\it Small
uncountable cardinals and topology,} Open problems
in topology (J. van Mill and G. Reed, eds.),
North--Holland, 1990, pp. 195-218.

\vfill\eject\end